\documentclass[a4paper,10pt,reqno]{amsart}

\usepackage[numbers]{natbib}
\usepackage{amsmath,amssymb}
\usepackage[foot]{amsaddr}
\usepackage{caption}
\usepackage{mathtools}
\usepackage[latin1]{inputenc}
\usepackage{mathrsfs}
\usepackage{enumerate}
\usepackage{textcmds}
\usepackage{graphicx,epstopdf}
\usepackage[caption=false]{subfig}
\usepackage[pdfusetitle]{hyperref}
\usepackage[capitalize]{cleveref}
\usepackage{booktabs}

\newtheorem{Theorem}{Theorem}
\newtheorem{Lemma}[Theorem]{Lemma}
\newtheorem{Remark}[Theorem]{Remark}
\newtheorem{Corollary}[Theorem]{Corollary}
\newenvironment{Proof}{\begin{proof}}{\end{proof}}

\newcommand{\D}[2]{\frac{\mathrm{d}#1}{\mathrm{d}#2}}
\newcommand{\T}{^{\operatorname{T}}}
\newcommand{\Osym}{\text{\usefont{OMS}{cmsy}{m}{n}O}}
\newcommand{\N}{\mathbb{N}}
\newcommand{\Z}{\mathbb{Z}}

\newcommand{\C}{\mathbb{C}}
\newcommand{\Cz}{\C^{2}}
\newcommand{\MzC}{\mathrm{M}_{2}(\C)}
\newcommand{\re}{\operatorname{Re}}
\newcommand{\spr}[2]{#2\T#1}

\title[On the connection coefficients for linear differential equations]{On the connection coefficients for linear\\differential systems with applications to the\\spheroidal and ellipsoidal wave equation}
\author{Harald Schmid}
\email{h.schmid@oth-aw.de}
\address{University of Applied Sciences Amberg-Weiden, Amberg, Germany}

\begin{document}

\begin{abstract}
This paper is concerned with the connection coefficients between the local fundamental solutions of a $2\times 2$ linear ordinary differential system with two neighboring regular singular points at $z=0$ and $z=1$. We derive an asymptotic formula for the connection coefficients which can be used for numerical calculations and, in particular, for determining the eigenvalues of some spectral problems arising in mathematical physics. As an application, new algorithms for computing the eigenvalues of the ellipsoidal wave equation and the spheroidal wave equation are presented.
\end{abstract}

\maketitle

\section{Introduction}

In \cite{SS:1980} and \cite{Schaefke:1980}, R. Sch\"{a}fke and D. Schmidt investigated the connection problem for a linear differential system
\begin{equation} \label{GenSys}
y'(z) = \left(\frac{1}{z}A + \frac{1}{z-1}B + G(z)\right)y(z)
\end{equation}
where $A$, $B$ are complex $n\times n$ matrices and $G(z)$ is a holomorphic matrix function on a disk $\{z\in\C:|z|<r\}$ with some $1<r\leq\infty$. This is the general form of a differential system with two neighboring regular singular points (i.e. singularities of the first kind) located at $0$ and $1$. Their main results are explicit asymptotic formulas for the connection coefficients between the fundamental sets of Floquet solutions. In the present paper we will focus on the case $n=2$. By an appropriate shift of the characteristic exponents and by multiplying the series coefficients with some suitably designed vector, we obtain an effective algorithm for the calculation of the connection coefficients. More precisely, in Section 2 we construct a sequence $\Theta_k$ that converges to the connection coefficient $\Theta$ to be determined, and we will give both the speed of convergence and an a posteriori estimate for the error of the approximation. Since many eigenvalue problems from physics and engineering can be reduced to the evaluation of the connection coefficients of an associated linear differential system, this approach also provides an algorithm for calculating the eigenvalues of some special functions. As an application of our result, a new algorithm for calculating the eigenvalues of the ellipsoidal wave equation is presented in Section 3. Moreover, in Section 4, we will increase the convergence speed of the algorithm introduced in \cite{Schmid:2023} for computing the eigenvalues of the spheroidal wave equation.

\section{Computation of the connection coefficients}

Let us consider the linear differential system \eqref{GenSys} for the case $n=2$, where $A,B\in\MzC$ are fixed matrices, and $G:\mathfrak{D}\longrightarrow\MzC$ is a holomorphic matrix function on the disk $\mathfrak{D}:=\{z\in\C:|z|<r\}$ with some $1<r\leq\infty$. We make the following assumptions:
\begin{enumerate}[(a)]
\item $a_0$ is an eigenvector of $A$ for the eigenvalue $\alpha_0$, and $A-\alpha_0-k$ is invertible for all positive integers $k$.
\item The matrix $B$ has distinct eigenvalues $\beta_1\neq\beta_2$, which are ordered in such a way that $\delta := \beta_2-\beta_1$ satisfies $\re(\delta)>-1$. Moreover, let $b_j$ be an eigenvector of $B$ for the eigenvalue $\beta_j$, $j\in\{1,2\}$.
\item The function $G(z)$ has the series expansion $G(z)=\sum_{k=0}^\infty z^k G_k$ on the unit disk $\mathfrak{D}_0:=\{z\in\C:|z|<1\}$, and $G(z)=\sum_{k=0}^\infty(1-z)^k\tilde G_k$ applies on the disk $\mathfrak{D}_1:=\{z\in\C:|z-1|<\min\{1,r-1\}\}$ centered at $z=1$.
\end{enumerate}
Condition (a) is fulfilled, for instance, if $\alpha_0$ is the maximum eigenvalue of $A$. According to (a), the system \eqref{GenSys} has a Floquet solution $y_0(z) = z^{\alpha_0}h_0(z)$ on $\mathfrak{G}_0:=\mathfrak{D}\setminus[1,r[$ with some holomorphic vector function $h_0:\mathfrak{G}_0\longrightarrow\Cz$ satisfying $h_0(0)=a_0$. In addition, there exists a fundamental matrix
\begin{equation} \label{FundMat}
Y(z) = H(z)
\begin{pmatrix} (1-z)^{\beta_1} & 0 \\[1ex] 0 & (1-z)^{\beta_2} \end{pmatrix}
\begin{pmatrix} 1 & 0 \\[1ex] q\log(1-z) & 1 \end{pmatrix}
\end{equation}
on $\mathfrak{G}_1:=\mathfrak{D}\setminus]{-r},0]$, where $H:\mathfrak{G}_1\longrightarrow\MzC$ is a holomorphic matrix function satisfying $H(1)e_1=b_1$ and $H(1)e_2=b_2$, and if $\beta_2-\beta_1$ is not an integer, then $q=0$ (cf. \cite[Theorem 5.6]{Wasow:1976}). Here and in the following, we denote by $e_1$, $e_2$ the standard basis of $\C^2$. Moreover, $z^{\alpha_0}$ is a uniquely determined holomorphic function on $\mathfrak{G}_1$ if we assume $-\pi<\arg z<\pi$, and $(1-z)^{\beta_j}$ is a holomorphic function on $\mathfrak{G}_0$ if $-\pi<\arg(1-z)<\pi$ is assumed. The fundamental matrix gives rise to the fundamental solutions
\begin{equation*}
y_1(z) := Y(z)e_1 = (1-z)^{\beta_1}h_1(z)+q(1-z)^{\beta_2}\log(1-z)h_2(z),\quad 
y_2(z) := Y(z)e_2 = (1-z)^{\beta_2}h_2(z)
\end{equation*}
on $\mathfrak{G}_1$ with some holomorphic vector functions $h_1,h_2:\mathfrak{G}_1\longrightarrow\MzC$. Now, on $\mathfrak{G}:=\mathfrak{G}_0\cap\mathfrak{G}_1$ we can write $y_0(z)$ as a linear combination
\begin{equation} \label{Connect}
y_0(z) = Y(z)c,\quad c = \begin{pmatrix} \Theta \\ \Omega \end{pmatrix}
\end{equation}
or equivalently $y_0(z) = \Theta y_1(z)+\Omega y_2(z)$ with $\Theta$, $\Omega$ called \emph{connection coefficients}. If $\Theta = 0$, then $y_0(z) = z^{\alpha_0}h_0(z)$ is a constant multiple of the Floquet solution $y_2(z)=(1-z)^{\beta_2}h_2(z)$. 

Many eigenvalue problems for linear Hamiltonian systems or second order differential equations (cf. \cite{Schmid:2023}) are reducible to the problem of finding one or more parameters such that certain Floquet solutions at different singular points can be smoothly matched. This in turn is equivalent to the problem of locating the zeros of one or more connection coefficients, for example those of $\Theta$.

\begin{Remark}
If the coefficient matrices $A=A(t)$, $B=B(t)$, $G(z)=G(z,t)$ depend analytically on some parameter $t\in\C$, where the eigenvalues $\alpha_0$ and $\beta_1$, $\beta_2$ itself are independent of $t$, then also $y_0(z)=y_0(z,t)$ and $Y(z)=Y(z,t)$ are analytical functions with respect to $t$ (see \cite[Lemma 6]{BSW:2005}). With an arbitrary but fixed $z_0\in\mathfrak{G}$ in \eqref{Connect}, we obtain $c(t)=Y(z_0,t)^{-1}y_0(z_0,t)$, so that $\Theta=\Theta(t)$ and $\Omega=\Omega(t)$ depend analytically on $t$ as well. In particular, the zeros of the entire functions $\Theta,\Omega:\C\longrightarrow\C$ form a discrete subset of $\C$.
\end{Remark}

In the following, we are looking for a method to compute $\Theta$ in an efficient manner. For this purpose we first apply the transformation
\begin{equation} \label{Shift}
y(z) = z^{\alpha_0}(1-z)^{\beta_1+1}\eta(z)
\end{equation}
which turns \eqref{GenSys} into the differential system
\begin{equation} \label{RegSing}
\eta'(z) = \left(\frac{1}{z}A_0 + \frac{1}{z-1}A_1 + G(z)\right)\eta(z)
\end{equation}
where $A_0 := A-\alpha_0$ and $A_1 := B-\beta_1-1$. It has the holomorphic solution
\begin{equation} \label{HolSol0}
\eta_0(z) := (1-z)^{-\beta_1-1}h_0(z) = \sum_{k=0}^\infty z^k d_k\quad\mbox{with}\quad d_0 = h_0(0) = a_0
\end{equation}
on $\mathfrak{D}_0$. Moreover, $\eta_j(z) := (1-z)^{-\beta_1-1}z^{-\alpha_0}y_j(z)$ for $j\in\{1,2\}$ are two fundamental solutions of \eqref{RegSing} on $\mathfrak{D}_1$, and \eqref{Connect} implies $\eta_0(z) = \Theta\eta_1(z) + \Omega\eta_2(z)$. As we will see below, we can compute $\Theta$ to be the limit of a sequence built directly from the series coefficients $d_k$. Thus, we will at first provide an algorithm that determines the vectors $d_k$.

\begin{Lemma} \label{thm:RecRel}
The coefficients $d_k$ of the holomorphic solution \eqref{HolSol0} can be computed by means of the two-step recurrence formula
\begin{equation} \label{Coeffs}
u_k := (A_0-k)^{-1}\left((A_1+1)d_{k-1} - \sum_{\ell=0}^{k-1} G_{k-1-\ell}u_{\ell}\right),\quad 
d_k := d_{k-1}+u_k\quad\mbox{for}\quad k=1,2,3,\ldots
\end{equation}
starting with $u_0 = d_0 := a_0$.
\end{Lemma}

\begin{Proof}
As $\eta_0(z)=\sum_{k=0}^\infty z^k d_k$ solves the system \eqref{RegSing}, we get
\begin{equation*}
A_0\eta_0(z) - z\eta_0'(z) = \left(z\big(A_0+A_1-G(z)\big) + z^2 G(z)\right)\eta_0(z) - z^2\eta_0'(z)
\end{equation*}
where the expression on the left hand side is $\sum_{k=0}^\infty z^k(A_0-k)d_k$. If we introduce $d_{-1} := 0$ to be the zero vector, then the right hand side becomes
\begin{align*}
& \sum_{k=0}^\infty z^{k+1}\left((A_0+A_1-k)d_k-\sum_{\ell=0}^k G_{k-\ell}d_{\ell}\right) + \sum_{k=0}^\infty z^{k+2}\left(\sum_{\ell=0}^k G_{k-\ell}d_{\ell}\right) \\
& = \sum_{k=0}^\infty z^k\left((A_0+A_1-k+1)d_{k-1} - \sum_{\ell=0}^{k-1} G_{k-1-\ell}(d_{\ell}-d_{\ell-1})\right)
\end{align*}
Comparing the coefficients yields the recurrence relation
\begin{align*}
(A_0-k)d_k = (A_0+A_1-k+1)d_{k-1} - \sum_{\ell=0}^{k-1}G_{k-1-\ell}(d_{\ell}-d_{\ell-1})
\end{align*}
If we define $u_k := d_k - d_{k-1}$ and $u_0 := d_0$, then we obtain
\begin{align*}
(A_0-k)u_k = (A_1+1)d_{k-1} - \sum_{\ell=0}^{k-1}G_{k-1-\ell}u_{\ell},\quad k=1,2,3,\ldots
\end{align*}
By means of (a), the matrix $A_0-k$ is invertible for all positive integers $k$. Multiplying above relation from the left by $(A_0-k)^{-1}$ implies \eqref{Coeffs}.
\end{Proof}

Using Lemma \ref{thm:RecRel}, we can also find the coefficients of the Floquet solution
\begin{equation} \label{FloqSol}
\eta_2(z) 
= (1-z)^{\beta_2-\beta_1-1}z^{-\alpha_0}h_2(z) 
= (1-z)^{\beta_2-\beta_1-1}\sum_{k=0}^\infty (1-z)^k\tilde d_k
\end{equation}
of \eqref{RegSing} starting with $\tilde d_0=h_2(1)=b_2$, since the transformation $\tilde\eta(z) = z^{1+\beta_1-\beta_2}\eta(1-z)$ changes \eqref{RegSing} to the system
\begin{equation} \label{RegSing1}
\tilde\eta'(z) = \left(\frac{1}{z}\tilde A_0 + \frac{1}{z-1}\tilde A_1 + \tilde G(z)\right)\tilde\eta(z)
\end{equation}
having the same structure as \eqref{RegSing}, where $\tilde A_0 := B-\beta_2$, $\tilde A_1 := A-\alpha_0$, and $\tilde G(z) := -G(1-z)=\sum_{k=0}^\infty z^k(-\tilde G_k)$. Hence, the coefficients $\tilde d_k$ of the holomorphic solution $\tilde\eta_0(z) := z^{1+\beta_1-\beta_2}\eta_2(1-z) = \sum_{k=0}^\infty z^k\tilde d_k$ can be computed with the recurrence formula \eqref{Coeffs} provided that we replace $A_0$, $A_1$, $G_k$, $d_k$ by $\tilde A_0$, $\tilde A_1$, $-\tilde G_k$, $\tilde d_k$ and that we start with $u_0 = \tilde d_0 := b_2$.

The following theorem shows how the series coefficients $d_k$, $\tilde d_k$ are related to the connection coefficient $\Theta$. Subsequently, we will use them to construct a sequence $\Theta_k$ that converges to $\Theta$ as $k\to\infty$.

\begin{Theorem} \label{thm:ConCoeff}
Let $d_k,\tilde d_k\in\Cz$ be the series coefficients of the holomorphic solution \eqref{HolSol0} and \eqref{FloqSol}, respectively. Moreover, if we set $\delta := \beta_2-\beta_1$, then for any fixed integer $n\geq 0$
\begin{equation} \label{LimCoeff}
d_k 
= \Theta b_1 + \frac{\omega}{\Gamma(k+1)}\sum_{\ell=0}^n(-1)^\ell\Gamma(\ell+\delta)\Gamma(k-\ell+1-\delta)\tilde d_\ell + \Osym(k^{-\re(\delta)-n-1})\quad\mbox{as}\quad k\to\infty
\end{equation}
with some constant $\omega$ which is independent of $k$ and $n$.
\end{Theorem}

\begin{Proof}
First, let us assume that $\beta_1$, $\beta_2$ do not differ by an integer, i.e. $\delta\not\in\Z$. In this case, $q=0$ in \eqref{FundMat}, so that
\begin{equation*}
\eta_1(z) = (1-z)^{-1}z^{-\alpha_0}h_1(z) = (1-z)^{-1}\sum_{k=0}^\infty (1-z)^k\tilde q_k
\end{equation*}
together with \eqref{FloqSol} forms a fundamental set of Floquet solutions on $\mathfrak{D}_1$, where $\tilde q_0=h_1(1)=b_1$ is an eigenvector of $B$ for the eigenvalue $\beta_1$. From \cite[Theorem 1.4]{SS:1980} with $\alpha=0$, $\alpha_1=-1$, $\alpha_2=\delta-1$ it follows that
\begin{equation} \label{AsymForm}
d_k 
= \Theta\sum_{\ell=0}^{n_1}\frac{\Gamma(k-\ell+1)}{\Gamma(k+1)\Gamma(-\ell+1)}\tilde q_\ell 
+ \Omega\sum_{\ell=0}^{n_2}\frac{\Gamma(k-\ell+1-\delta)}{\Gamma(k+1)\Gamma(-\ell+1-\delta)}\tilde d_\ell
+ \Osym(k^{-\nu-1})\quad\mbox{as}\quad k\to\infty
\end{equation}
where $n_1,n_2\geq 0$ are arbitrary integers and $\nu := \min\{n_1,\re(\delta)+n_2\}$. Since $\frac{1}{\Gamma(-\ell+1)}=0$ for any integer $\ell>0$, we get
\begin{equation*}
d_k = \Theta\frac{\Gamma(k+1)}{\Gamma(k+1)\Gamma(1)}b_1 
+ \Omega\sum_{\ell=0}^{n_2}\frac{\Gamma(k-\ell+1-\delta)}{\Gamma(k+1)\Gamma(-\ell+1-\delta)}\tilde d_\ell + \Osym(k^{-\nu-1})
\end{equation*}
being independent of $n_1$. If we choose $n_1$ sufficiently large and if we set $n_2=n$, then $\nu = \re(\delta)+n$. Moreover, as $\frac{1}{\Gamma(-\ell+1-\delta)} = (-1)^\ell\frac{\sin(\pi\delta)}{\pi}\Gamma(\ell+\delta)$, we obtain the asymptotic formula \eqref{LimCoeff} with $\omega:=\frac{1}{\pi}\Omega\sin(\pi\delta)$.

It should be noted that the previously applied transformation \eqref{Shift} has significant impact on the asymptotic formula \eqref{AsymForm} for $d_k$. By shifting the Floquet exponents to $\alpha_1=-1$ and $\alpha=0$ in \cite[Theorem 1.4]{SS:1980}, only one term $\Theta b_1$ is left from the first sum, and the terms in the second sum are easier to handle. 

Now, let us consider the case where $\delta=m$ is a positive integer (because of assumption (b) we can omit $\delta=0$). Here, \eqref{FundMat} takes the form
\begin{align*}
Y(z) & = H(z)
\begin{pmatrix} (1-z)^{\beta_1} & 0 \\[1ex] 0 & (1-z)^{\beta_1+m} \end{pmatrix}
\begin{pmatrix} 1 & 0 \\[1ex] q\log(1-z) & 1 \end{pmatrix}
\end{align*} 
and by means of \eqref{Shift},
\begin{equation*}
\hat Y(z) := (1-z)^{-\beta_1-1}z^{-\alpha_0}Y(z) 
= (1-z)^{-1}\hat H(z)\begin{pmatrix} 1 & 0 \\[1ex] q\log(1-z) & 1 \end{pmatrix}
\end{equation*}
is a fundamental matrix of \eqref{RegSing}, where
\begin{equation*}
\hat H(z) := z^{-\alpha_0}H_1(z)\begin{pmatrix} 1 & 0 \\[1ex] 0 & (1-z)^m \end{pmatrix} 
= \sum_{k=0}^\infty (1-z)^k D_k
\end{equation*}
is a holomorphic matrix function on $\mathfrak{D}_1$. Its coefficients $D_k\in\MzC$ have the form
\begin{equation} \label{AsymExp}
D_k = \begin{pmatrix} \ast & 0 \\[1ex] \ast & 0 \end{pmatrix}\quad\mbox{for}\quad k=0,\ldots,m-1
\end{equation}
where $D_0 e_1 = H(1)e_1 = b_1$ and $D_m e_2 = H(1)e_2 = b_2$. Moreover, since
\begin{equation*}
(1-z)^{-1}\begin{pmatrix} 1 & 0 \\[1ex] q\log(1-z) & 1 \end{pmatrix} = (1-z)^Q
\quad\mbox{with}\quad Q := \begin{pmatrix} -1 & 0 \\[1ex] q & -1 \end{pmatrix}
\end{equation*}
we obtain $\hat Y(z) = \sum_{k=0}^\infty (1-z)^k D_k(1-z)^Q$. Applying \cite[Theorem 1.1]{Schaefke:1980} with $\alpha=0$ and $\gamma_{-}=-1$, it follows that
\begin{equation} \label{dApprox}
d_k = \sum_{\ell=0}^N\frac{1}{\Gamma(k+1)}D_\ell\frac{1}{\Gamma}(-\ell-Q)\Gamma(k-\ell-Q)c + \Osym(k^{-N-1+\varepsilon})\quad\mbox{as}\quad k\to\infty
\end{equation}
for arbitrary $N\in\N$ and $\varepsilon>0$. In particular, for any integers $\ell>0$ and $k\geq\ell$ we get
\begin{align*}
\frac{1}{\Gamma}(-\ell-Q) 
& = \begin{pmatrix} \frac{1}{\Gamma}(1-\ell) & 0 \\[1ex] -q\big(\frac{1}{\Gamma}\big)'(1-\ell) & \frac{1}{\Gamma}(1-\ell) \end{pmatrix} 
  = \begin{pmatrix} 0 & 0 \\[1ex] q (-1)^{\ell}(\ell-1)! & 0 \end{pmatrix} \\
\Gamma(k-\ell-Q) 
& = \begin{pmatrix} \frac{1}{\Gamma}(k+1-\ell) & 0 \\[1ex] -q\big(\frac{1}{\Gamma}\big)'(k+1-\ell) & \frac{1}{\Gamma}(k+1-\ell) \end{pmatrix}^{-1}
  = (k-\ell)!\begin{pmatrix} 1 & 0 \\[1ex] -q\psi(k+1-\ell) & 1 \end{pmatrix}
\end{align*}
according to \cite[Theorem A.2.(ii) and A.3]{Schaefke:1980}, where $\psi(z):=\D{}{z}\ln\Gamma(z)$ denotes the digamma function. Hence,
\begin{equation} \label{Term1}
D_{\ell}\frac{1}{\Gamma}(-\ell-Q)\Gamma(k-\ell-Q) 
= q(-1)^{\ell}(\ell-1)!(k-\ell)! D_{\ell}\begin{pmatrix} 0 & 0 \\[1ex] 1 & 0 \end{pmatrix}
\end{equation}
for all integers $k\geq\ell>0$. On the other hand, if $\ell=0$, then we get
\begin{align*}
\frac{1}{\Gamma}(-Q) 
& = \begin{pmatrix} \frac{1}{\Gamma}(1) & 0 \\[1ex] -q\big(\frac{1}{\Gamma}\big)'(1) & \frac{1}{\Gamma}(1) \end{pmatrix} = \begin{pmatrix} 1 & 0 \\[1ex] -q\gamma & 1 \end{pmatrix},\quad
\Gamma(k-Q) =  k!\begin{pmatrix} 1 & 0 \\[1ex] -q\psi(k+1) & 1 \end{pmatrix}
\end{align*}
where $\gamma=0.577215\ldots$ is Euler's constant. From \eqref{AsymExp} it follows that
\begin{equation} \label{OTerms}
D_0\frac{1}{\Gamma}(-Q)\Gamma(k-Q) = k! D_0,\quad
D_{\ell}\frac{1}{\Gamma}(-\ell-Q)\Gamma(k-\ell-Q)=O\quad\mbox{for}\quad\ell=1,\ldots,m-1
\end{equation}
Moreover, if we choose $N=n+m+1$ and $\varepsilon=1$, then \eqref{dApprox}, \eqref{Term1} and \eqref{OTerms} yield
\begin{align*}
d_k 
& = D_0 c + \sum_{\ell=m}^{m+n+1}q\frac{(-1)^{\ell}(\ell-1)!(k-\ell)!}{\Gamma(k+1)}D_{\ell}\begin{pmatrix} 0 & 0 \\[1ex] 1 & 0 \end{pmatrix}c + \Osym(k^{-m-n-1}) \\
& = \Theta b_1 + q\Theta\sum_{\ell=0}^{n+1}\frac{(-1)^{\ell+m}(\ell+m-1)!(k-\ell-m)!}{\Gamma(k+1)}D_{\ell+m}e_2 + \Osym(k^{-m-n-1})
\end{align*}
Note that the Floquet solution \eqref{FloqSol} can be written in the form
\begin{equation*}
\eta_2(z) = \hat Y(z)e_2 
= (1-z)^{-1}\sum_{k=m}^\infty (1-z)^k D_k e_2 
= (1-z)^{m-1}\sum_{k=0}^\infty (1-z)^k D_{k+m} e_2
\end{equation*}
and hence we get $D_{k+m}e_2=\tilde d_k$ for all $k$. Since $m=\delta$, we have $(\ell+m-1)!=\Gamma(\ell+\delta)$ and $(k-\ell-m)!=\Gamma(k-\ell+1-\delta)$. Finally, if we set $\omega:=(-1)^m q\Theta$, then
\begin{equation*}
d_k 
= \Theta b_1 + \frac{\omega}{\Gamma(k+1)}\sum_{\ell=0}^{n+1}(-1)^\ell\Gamma(\ell+\delta)\Gamma(k-\ell+1-\delta)\tilde d_\ell + \Osym(k^{-\delta-n-1})
\end{equation*}
From $\frac{\Gamma(k-n-\delta)}{\Gamma(k+1)} = k^{-\delta-n-1}(1+\Osym(k^{-1}))$, cf. \cite[Section 1.1]{MOS:1966}, it follows that the last term in the sum for $\ell=n+1$ as well as the remainder asymptotically behave like $\Osym(k^{-\delta-n-1})$. Thus, \eqref{LimCoeff} holds even if $\delta$ is a positive integer, and this completes the proof of Theorem \ref{thm:ConCoeff}.
\end{Proof}

In a next step, for a given integer $n\geq 0$ and for all integers $k>\re(\delta)+n-1$ we introduce the vectors
\begin{equation} \label{Sum}
p_k 
:= b_2 + \sum_{\ell=1}^n \left(\prod_{m=0}^{\ell-1}\frac{m+\delta}{m+\delta-k}\right)\tilde d_\ell
 = b_2 + \Osym\big(\tfrac{1}{k}\big)
\end{equation}
Since $b_1$ and $b_2$ are linearly independent, there exists a constant $k_1$ such that $b_1$ and $p_k$ are also linearly independent for all $k\geq k_1$, and we can define the vectors
\begin{equation} \label{Ortho}
\vartheta_k := \frac{1}{\spr{Jp_k}{b_1}}Jp_k,\quad\mbox{where}\quad 
J := \begin{pmatrix} 0 & 1 \\[1ex] -1 & 0 \end{pmatrix}
\end{equation}
which satisfy $\spr{b_1}{\vartheta_k} = 1$ and $\spr{p_k}{\vartheta_k} = 0$ for all $k\geq k_1$ (note that $\spr{Jp_k}{b_1}=\det(b_1,p_k)\neq 0$).

\begin{Theorem} \label{thm:Theta}
Let $d_k,\tilde d_k\in\Cz$ be the series coefficients of the holomorphic solutions \eqref{HolSol0} and \eqref{FloqSol}, respectively. In addition, let $n$ be a given nonnegative integer. Taking the vectors $\vartheta_k$ for $k\geq k_1$ defined by \eqref{Ortho}, we get
\begin{equation} \label{Theta}
\Theta_k := \spr{d_k}{\vartheta_k} = \Theta + \Osym(k^{-\re(\delta)-n-1})\quad\mbox{as}\quad k\to\infty
\end{equation}
Furthermore, the limit
\begin{equation*}
\tau := \lim_{k\to\infty} k^{\delta+n+2}(\Theta_k-\Theta_{k-1})
\end{equation*}
exists, and for an arbitrary $\varepsilon>0$ there is an integer $k_2$ such that
\begin{equation} \label{ConRate}
|\Theta - \Theta_k| \leq (1+\varepsilon)\frac{k\,|\Theta_k-\Theta_{k-1}|}{\re(\delta)+n+1}
\quad\mbox{for all}\quad k\geq k_2
\end{equation}
\end{Theorem}

\begin{Proof}
Substituting $(-1)^\ell\Gamma(\ell+\delta)\Gamma(k-\ell+1-\delta) = \Gamma(\delta)\Gamma(k+1-\delta)\prod_{m=0}^{\ell-1}\frac{m+\delta}{m+\delta-k}$ in \eqref{LimCoeff} gives
\begin{equation} \label{dCoeff}
d_k = \Theta b_1 + \omega\frac{\Gamma(\delta)\Gamma(k+1-\delta)}{\Gamma(k+1)}p_k + \Osym(k^{-\re(\delta)-n-1})\quad\mbox{as}\quad k\to\infty
\end{equation}
Since $p_k = b_2 + \tfrac{\delta}{\delta-k}\tilde d_1 + \Osym(k^{-2})$, the associated vector \eqref{Ortho} takes the form $\vartheta_k = \zeta_0 + \tfrac{1}{k}\zeta_1 + \Osym(k^{-2})$, where
\begin{equation*}
\zeta_0 := \frac{1}{\spr{Jb_2}{b_1}}Jb_2
\end{equation*}
and the vector $\zeta_1$ is independent of $k\geq k_1$. Multiplying \eqref{dCoeff} from the left by $\vartheta_k\T$ immediately yields \eqref{Theta}. If we additionally take into account the series coefficient $\tilde d_{n+1},\ldots,\tilde d_N$ in \eqref{LimCoeff} for some integer $N\geq n+1$, then
\begin{equation*}
\Theta_k = \Theta + 
\omega\sum_{\ell=n+1}^N(-1)^\ell\frac{\Gamma(\ell+\delta)\Gamma(k-\ell+1-\delta)}{\Gamma(k+1)}\spr{\tilde d_\ell}{\vartheta_k} + \Osym(k^{-\re(\delta)-N-1})
\end{equation*}
and hence
\begin{equation*}
\Theta_k-\Theta_{k-1} = \omega\sum_{\ell=n+1}^N (-1)^\ell\Gamma(\ell+\delta)\spr{\tilde d_{\ell}}{\Delta_{k,\ell}} + \Osym(k^{-\re(\delta)-N-1})
\end{equation*}
where
\begin{equation*}
\Delta_{k,\ell} 
:= \frac{\Gamma(k-\ell+1-\delta)}{\Gamma(k+1)}\vartheta_k - \frac{\Gamma(k-\ell-\delta)}{\Gamma(k)}\vartheta_{k-1} 
 = \frac{\Gamma(k-\ell-\delta)}{\Gamma(k)}\left(-\frac{\ell+\delta}{k}\vartheta_k + \vartheta_k - \vartheta_{k-1}\right)
\end{equation*}
Since $\vartheta_k - \vartheta_{k-1} = \Osym(k^{-2})$ and $\frac{\Gamma(k-\ell-\delta)}{\Gamma(k)} = k^{-\delta-\ell}(1+\Osym(k^{-1}))$, cf. \cite[Section 1.1]{MOS:1966}, we get
\begin{equation*}
\Gamma(\ell+\delta)\Delta_{k,\ell} = k^{-\delta-\ell-1}\big(-\Gamma(\ell+\delta+1)\zeta_0+\Osym\big(\tfrac{1}{k}\big)\big)
\end{equation*}
In particular, if we choose $N = n+2$, then
\begin{align*}
\Theta_k-\Theta_{k-1} 
& = \omega\sum_{\ell=n+1}^{n+2} (-1)^{\ell} k^{-\delta-\ell-1}\big(-\Gamma(\ell+\delta+1)\spr{\tilde d_\ell}{\zeta_0} + \Osym\big(\tfrac{1}{k}\big)\big) + \Osym(k^{-\re(\delta)-n-3}) \\
& = k^{-\delta-n-2}\left(\omega\,(-1)^{n}\Gamma(n+2+\delta)\spr{\tilde d_{n+1}}{\zeta_0} + \Osym\big(\tfrac{1}{k}\big)\right)\quad\mbox{as}\quad k\to\infty
\end{align*}
and therefore the limit
\begin{equation*}
\tau = \lim_{k\to\infty} k^{\delta+n+2}(\Theta_k-\Theta_{k-1}) = \omega\,(-1)^{n}\Gamma(n+2+\delta)\spr{\tilde d_{n+1}}{\zeta_0}
\end{equation*}
exists. Now, if we set $\varepsilon_0:=\frac{\varepsilon}{2+\varepsilon}$ for a given $\varepsilon>0$, then $0<\varepsilon_0<1$, and we can find some integer $k_2\geq k_1$ such that 
$|k^{\delta+n+2}(\Theta_k-\Theta_{k-1})-\tau|\leq\varepsilon_0|\tau|$ holds for all $k\geq k_2$. This means that $|\Theta_k-\Theta_{k-1}|\leq(1+\varepsilon_0)|\tau|k^{-\re(\delta)-n-2}$ is satisfied for all $k\geq k_2$, and we obtain
\begin{equation*}
|\Theta - \Theta_k| 
= \left|\sum_{\kappa=k+1}^\infty \Theta_{\kappa}-\Theta_{\kappa-1}\right|
\leq\sum_{\kappa=k+1}^\infty\frac{(1+\varepsilon_0)|\tau|}{\kappa^{\re(\delta)+n+2}}
\leq\int_k^\infty\frac{(1+\varepsilon_0)|\tau|}{\kappa^{\re(\delta)+n+2}}\,\mathrm{d}\kappa
= \frac{(1+\varepsilon_0)|\tau|}{\re(\delta)+n+1}k^{-\re(\delta)-n-1}
\end{equation*}
On the other hand, $|\tau|\leq\varepsilon_0|\tau|+|k^{\delta+n+2}(\Theta_k-\Theta_{k-1})|$ implies $|\tau|\leq\frac{1}{1-\varepsilon_0}k^{\re(\delta)+n+2}|\Theta_k-\Theta_{k-1}|$, and therefore
\begin{equation*}
|\Theta - \Theta_k| \leq \frac{1+\varepsilon_0}{1-\varepsilon_0}\frac{k\,|\Theta_k-\Theta_{k-1}|}{\re(\delta)+n+1}\quad\mbox{for all}\quad k\geq k_2
\end{equation*}
Finally, using $\frac{1+\varepsilon_0}{1-\varepsilon_0}=1+\varepsilon$ yields the estimate \eqref{ConRate}.
\end{Proof}

Since $\re(\delta)>-1$ by assumption (b), we have $-\re(\delta)-n-1<0$, and from \eqref{Theta} it follows that $\lim_{k\to\infty}\Theta_k = \Theta$ for any choice of $n\geq 0$. Furthermore, \eqref{ConRate} gives an a posteriori estimate for the error $|\Theta-\Theta_k|$ provided that $k^{\delta+2+n}(\Theta_k-\Theta_{k-1})$ is reasonably close to its limit $\tau$. 

The vector $\vartheta_k$ used in Theorem \ref{thm:Theta} results from a finite sum \eqref{Sum} with $n+1$ terms. The more terms we involve to build $\vartheta_k$, the faster the sequence $\Theta_k$ converges to $\Theta$ as $k\to\infty$. More precisely, the length $n+1$ of the sum is related to the speed of convergence via $\Osym(k^{-\re(\delta)-n-1})$. In the case $n=0$, we have $p_k = b_2$ for all $k$, and thus $\vartheta_k = \zeta_0 = (\spr{Jb_2}{b_1})^{-1}Jb_2$ needs to be calculated only once, but on the other hand for $n=0$ we get the lowest speed of convergence $\Osym(k^{-\re(\delta)-1})$.

\section{Application to the ellipsoidal wave equation}

In this section we study the angular ellipsoidal wave equation (or Lam\'{e} wave equation) in algebraic form
\begin{equation} \label{AEWE}
z(z-1)(z-c)w''(z) + \tfrac{1}{2}(3z^2-2(1+c)z+c)w'(z)+(\lambda+\mu z + \gamma z^2)w(z)=0
\end{equation} 
It is obtained by separating the Helmholtz equation $\nabla^2\psi + k^2\psi = 0$ in ellipsoidal coordinates (see \cite{ATZ:1983}, for example). In this context, $c>1$ and $\gamma$ are real numbers which are related to $k^2$ and to the geometry of the ellipsoidal coordinates, while $\lambda$ and $\mu$ are two numbers resulting from the separation of $\psi$; they are called the \emph{eigenvalue parameters} of \eqref{AEWE}. 

In contrast to other special functions from mathematical physics, ellipsoidal wave functions have not yet been extensively studied (except for the special case $\gamma=0$, which is the so-called Lam\'{e} equation). Nevertheless, in \cite{Fedoryuk:1989a} and \cite{Fedoryuk:1989b} some basic properties of the ellipsoidal wave equation are established: It is shown that the spectrum is real, and estimates for the eigenvalues $(\lambda,\mu)$ as well as orthogonality relations for the associated Lam\'{e} wave functions are given. The number of papers dealing with the numerical computation of these eigenvalues is likewise rather small. One such approach, which is due to Arscott et\,al. (see \cite{ATZ:1983} and \cite{WL:2005}), results in a third-order difference equation whose minimal solution is approximated by backward recurrence, with a nonlinear system of equations to be finally solved. The computational difficulties arising from this method are also discussed in \cite{ATZ:1983}. A completely different approach is employed by Abramov et\,al. in \cite{ADKL:1989}, where the eigenvalues are determined using an auxiliary differential equation for an associated phase function. 

By using the results found in the last section, we are now going to derive an efficient method for the numerical computation of the eigenvalues and eigenfunctions of the ellipsoidal wave functions \eqref{AEWE}. To this end, we consider the $2\times 2$ linear differential system
\begin{equation} \label{EllSys}
y'(z) = \left(\frac{1}{z}\begin{pmatrix} -\frac{1}{2} & a_{12} \\[1ex] 0 & 0 \end{pmatrix} 
+ \frac{1}{z-1}\begin{pmatrix} -\frac{1}{2} & b_{12} \\[1ex] 0 & 0 \end{pmatrix} 
+ \frac{1}{z-c}\begin{pmatrix} -\frac{1}{2} & r_{12} \\[1ex] 0 & 0 \end{pmatrix} 
- \frac{1}{c}\begin{pmatrix} 0 & 0 \\[1ex] 1 & 0 \end{pmatrix}\right)y(z)
\end{equation}
which is associated to the ellipsoidal wave equation in the following way: If $w(z)$ is a solution of \eqref{AEWE} and if we set
\begin{equation} \label{Param}
y(z) := \begin{pmatrix} -cw'(z) \\[0.5ex] w(z) \end{pmatrix},\quad
a_{12} := \lambda,\quad
b_{12} := \frac{c(\lambda+\mu+\gamma)}{1-c},\quad
r_{12} := \frac{\lambda+c\mu+c^2\gamma}{c-1}
\end{equation}
then $y(z)$ is a solution of \eqref{EllSys}. Conversely, if $y(z)$ solves this system for some parameters $a_{12},b_{12},r_{12}\in\C$, then we get a solution of the ellipsoidal wave equation \eqref{AEWE} by means of
\begin{equation*}
w(z) := \spr{y(z)}{e_2},\quad
\lambda := a_{12},\quad
\mu := -\frac{c(a_{12}+b_{12})+a_{12}+r_{12}}{c},\quad
\gamma := \frac{a_{12}+b_{12}+r_{12}}{c}
\end{equation*}
The differential system \eqref{EllSys} has three regular-singular points at $z_0=0$, $z_1=1$ and $z_2=c$, all with characteristic exponents $0$ or $-\frac{1}{2}$, from which we get six Floquet solutions
\begin{equation*}
y_{i,j}(z)=(z-z_i)^{-j/2}\sum_{k=0}^\infty(z-z_i)^k d_{i,j,k}
\end{equation*}
having the form
\begin{equation} \label{FundSet}
y_{i,0}(z) = \begin{pmatrix} \ast \\[1ex] 1 \end{pmatrix}  + \Osym(z-z_i),\quad
y_{i,1}(z) = (z-z_i)^{-1/2}\left(\begin{pmatrix} \frac{1}{2} \\[1ex] 0 \end{pmatrix} 
+ (z-z_i)\begin{pmatrix} \ast \\[1ex] 1 \end{pmatrix} + \Osym\big((z-z_i)^2\big)\right)
\end{equation}
for $i\in\{0,1,2\}$ and $j\in\{0,1\}$ on $\mathfrak{D}_i\subset\C$, where 
\begin{equation*}
\mathfrak{D}_0 := \{z\in\C:|z|<1\},\quad
\mathfrak{D}_1 := \{z\in\C:|z-1|<\min\{1,c-1\}\},\quad
\mathfrak{D}_2 := \{z\in\C:|z-c|<c-1\}
\end{equation*}
Thus, the ellipsoidal wave equation \eqref{AEWE} has six solutions $w_{i,j}(z) := \spr{y_{i,j}(z)}{e_2}$ for $i\in\{0,1,2\}$ and $j\in\{0,1\}$, which behave like $w_{i,j}(z) = (z-z_i)^{j/2}(1+\Osym(z-z_i))$ as $z\to z_i$. If \eqref{AEWE} has a nontrivial solution of the form 
\begin{equation*}
w(z) = z^{\rho/2}(z-1)^{\sigma/2}(z-c)^{\tau/2}f(z)
\end{equation*}
for given $\rho,\sigma,\tau\in\{0,1\}$ with some entire function $f:\C\longrightarrow\C$, then $w(z)$ is called \emph{ellipsoidal wave function} (cf. \cite[Sec. 3]{ATZ:1983} or \cite[\S 2]{Fedoryuk:1989b}). In this case, the functions $w_{0,\rho}$, $w_{1,\sigma}$ as well as $w_{1,\sigma}$, $w_{2,\tau}$ need to be linearly dependent on $\mathfrak{D}_0\cap\mathfrak{D}_1$ and $\mathfrak{D}_1\cap\mathfrak{D}_2$, respectively. Hence, the eigenvalue problem for the ellipsoidal wave equation consists of eight distinct spectral problems, and in contrast to classical problems for Sturm-Liouville equations, they contain two spectral parameters $\lambda$ and $\mu$. These values have to be determined such that $w_{0,\rho} = \Omega w_{1,\sigma}$ and $w_{2,\tau} = \hat\Omega w_{1,\sigma}$ or equivalently $y_{0,\rho} = \Omega y_{1,\sigma}$ and $y_{2,\tau} = \hat\Omega y_{1,\sigma}$ hold with some constants $\Omega,\hat\Omega\in\C$. 

Nevertheless, it suffices to focus on the connection between the Floquet solutions at the regular singular points $0$ and $1$ due to the following reason: The transformation 
\begin{equation} \label{Trafo}
\hat y(\hat z) := \begin{pmatrix} -1 & 0 \\[1ex] 0 & 1 \end{pmatrix}y(z),\quad z=c+(1-c)\hat z
\end{equation}
turns \eqref{EllSys} into the system
\begin{equation} \label{AEWS}
\hat y'(\hat z) = \left(
  \frac{1}{\hat z}\begin{pmatrix} -\frac{1}{2} & \hat a_{12} \\[1ex] 0 & 0 \end{pmatrix} 
+ \frac{1}{\hat z-1}\begin{pmatrix} -\frac{1}{2} & \hat b_{12} \\[1ex] 0 & 0 \end{pmatrix} 
+ \frac{1}{\hat z-\hat c}\begin{pmatrix} -\frac{1}{2} & \hat r_{12} \\[1ex] 0 & 0 \end{pmatrix}
- \frac{1}{\hat c}\begin{pmatrix} 0 & 0 \\[1ex] 1 & 0 \end{pmatrix}\right)\hat y(\hat z)
\end{equation}
having the same structure as \eqref{EllSys} with
\begin{equation} \label{NewPar}
\hat a_{12} := -r_{12},\quad
\hat b_{12} := -b_{12},\quad
\hat r_{12} := -a_{12},\quad
\hat c := \frac{c}{c-1} > 1
\end{equation}
where $z=1$ and $z=c$ are mapped to the regular singular points $\hat z=1$ and $\hat z=0$, respectively. Now, if $\hat y_{i,j}$ are the associated Floquet solutions of \eqref{AEWS} having the Form \eqref{FundSet}, then $y_{2,\tau} = \hat\Omega y_{1,\sigma}$ is equivalent to $\hat y_{0,\tau} = \hat\Omega\hat y_{1,\sigma}$. Since in general $y_{0,\rho} = \Theta y_{1,1-\sigma} + \Omega y_{1,\sigma}$ and $\hat y_{0,\tau} = \hat\Theta\hat y_{1,1-\sigma} + \hat\Omega\hat y_{1, \sigma}$ holds, we need to find values $\lambda$ and $\mu$ such that $\Theta(\lambda,\mu)=0$ and $\hat\Theta(\lambda,\mu)=0$ are satisfied simultaneously. In the following we deal with the computation of these connection coefficients. We can restrict our attention to $\Theta$, since for the calculation of $\hat\Theta$ we only have to exchange $\rho$ by $\tau$ and the entries in \eqref{EllSys} by \eqref{NewPar}.

The differential system \eqref{EllSys} has the form
\begin{equation*}
y'(z) = \left(\frac{1}{z}A + \frac{1}{z-1}B + G(z)\right)y(z),\quad
G(z) := \frac{1}{z-c}R-\frac{1}{c}S = \sum_{k=0}^\infty z^k G_k
\end{equation*}
where $G_0 := -\frac{1}{c}(R+S)$, $G_k := -\frac{1}{c^{k+1}}R$ for $k>0$, and
\begin{equation*}
A := \begin{pmatrix} -\frac{1}{2} & a_{12} \\[1ex] 0 & 0 \end{pmatrix},\quad
B := \begin{pmatrix} -\frac{1}{2} & b_{12} \\[1ex] 0 & 0 \end{pmatrix},\quad
R := \begin{pmatrix} -\frac{1}{2} & r_{12} \\[1ex] 0 & 0 \end{pmatrix},\quad
S := \begin{pmatrix} 0 & 0 \\[1ex] 1 & 0 \end{pmatrix}
\end{equation*}
Therefore, it is a special case of \eqref{GenSys}. In order to smoothly match the Floquet solutions $y_{0,\rho}$ and $y_{1,\sigma}$ for some $\rho,\sigma\in\{0,1\}$, we have to consider the eigenvalue $\alpha_0=-\rho/2$ of $A$ and the eigenvalues $\beta_1=(\sigma-1)/2$, $\beta_2=-\sigma/2$ of $B$ with corresponding eigenvectors
\begin{equation*}
a_0 := \begin{pmatrix} 2a_{12}(1-\rho)+\frac{\rho}{2} \\[1ex] 1-\rho \end{pmatrix},\quad
b_1 := \begin{pmatrix} 2b_{12}\sigma+\frac{1-\sigma}{2} \\[1ex] \sigma \end{pmatrix},\quad
b_2 := \begin{pmatrix} 2b_{12}(1-\sigma)+\frac{\sigma}{2} \\[1ex] 1-\sigma \end{pmatrix}
\end{equation*}
Since $\delta := \beta_2-\beta_1 = 1/2-\sigma>-1$, the assumptions (a) -- (c) in Section 2 are fulfilled. Moreover, if we set $A_0 := A-\alpha_0 = A+\rho/2$ and $A_1 := B-\beta_1-1 = B-(\sigma+1)/2$, then the recurrence relation \eqref{Coeffs} for computing the coefficients of the holomorphic solution $\eta_0(z)=\sum_{k=0}^\infty z^k d_k$ of \eqref{RegSing} reads
\begin{equation}\begin{split} \label{RecRel0}
u_k & := \big(A+\tfrac{\rho}{2}-k\big)^{-1}\left(\big(B-\tfrac{\sigma}{2}+\tfrac{1}{2}\big)d_{k-1} + \tfrac{1}{c}Su_{k-1} + \tfrac{1}{c}Rs_{k-1}\right) \\
d_k & := d_{k-1}+u_k,\quad s_k := \tfrac{1}{c}s_{k-1}+u_k\quad\mbox{for}\quad k=1,2,3,\ldots
\end{split}\end{equation}
starting with $u_0 = s_0 = d_0 := a_0$. Finally, we need to calculate the coefficients of the Floquet solution $\tilde\eta_0(z) = \sum_{k=0}^\infty z^k\tilde d_k$ of \eqref{RegSing1}, where $\tilde A_0 = B-\beta_2 = B+\sigma/2$, $\tilde A_1=A-\alpha_0=A+\rho/2$, and 
\begin{equation*}
\tilde G(z) = -G(1-z) = \frac{1}{c}S-\frac{1}{1-c-z}R 
= \frac{1}{c}S-\sum_{k=0}^\infty\frac{z^k}{(1-c)^{k+1}}R
\end{equation*}
With a similar reasoning as above, we obtain the coefficients from the recurrence formula
\begin{equation}\begin{split} \label{RecRel1}
\tilde u_k & := \big(B+\tfrac{\sigma}{2}-k\big)^{-1}\left(\big(A+\tfrac{\rho}{2}+1\big)\tilde d_{k-1} - \tfrac{1}{c}S\tilde u_{k-1} + \tfrac{1}{1-c}R\tilde s_{k-1}\right) \\
\tilde d_k & := \tilde d_{k-1}+\tilde u_k,\quad\tilde s_k := \tfrac{1}{1-c}\tilde s_{k-1}+\tilde u_k\quad\mbox{for}\quad k=1,2,3,\ldots
\end{split}\end{equation}
where $\tilde u_0 = \tilde s_0 = \tilde d_0 := b_2$. The following assertion follows directly from Theorem \ref{thm:Theta} with $\delta = 1/2-\sigma$:

\begin{Theorem} \label{thm:AEWF}
For fixed $\rho,\sigma\in\{0,1\}$ and for the parameters $a_{12}$, $b_{12}$, $r_{12}$ given by \eqref{Param}, let $d_k$ and $\tilde d_k$ for $k=0,1,2,\ldots$ be the vectors defined by the recurrence relations \eqref{RecRel0} and \eqref{RecRel1}, respectively. Moreover, for an arbitrary $n\geq 0$, let
\begin{equation*}
p_k := b_2 + \sum_{\ell=1}^n \left(\prod_{m=0}^{\ell-1}\frac{m+\frac{1}{2}-\sigma}{m+\frac{1}{2}-\sigma-k}\right)\tilde d_\ell,\quad 
\vartheta_k := \frac{1}{b_1\T Jp_k}Jp_k,\quad 
J := \begin{pmatrix} 0 & 1 \\[1ex] -1 & 0 \end{pmatrix}
\end{equation*}
where $b_1\T Jp_k\neq 0$ for sufficiently large $k$, say for $k\geq k_1$. If we define $\Theta_k := \spr{d_k}{\vartheta_k}$ for all $k\geq k_1$, then $\lim_{k\to\infty}\Theta_k = \Theta$, and the convergence being $\Osym(k^{\sigma-n-3/2})$.
\end{Theorem}

For given parameters $\lambda$ and $\mu$, the value $\Theta(\lambda,\mu)$ can be calculated with arbitrarily high precision using the sequence $\Theta_k$ in Theorem \ref{thm:AEWF}, where additionally \eqref{ConRate} yields the a posteriori estimate 
\begin{equation*}
|\Theta - \Theta_k| \leq (1+\varepsilon)\frac{k\,|\Theta_k-\Theta_{k-1}|}{n+1.5-\sigma}
\end{equation*}
for sufficiently large $k\geq k_2$ and arbitrary small $\varepsilon>0$. Replacing the values $a_{12}$, $b_{12}$, $r_{12}$, $c$ by \eqref{NewPar} and $\rho$ by $\tau$, we can apply Theorem \ref{thm:AEWF} to determine the connection coefficient $\hat\Theta(\lambda,\mu)$ as well, provided that we generate corresponding sequences $\hat d_k$ and $\hat\Theta_k$ for the system \eqref{AEWS}. When computing $\Theta$ or $\hat\Theta$, the value $n$ should be chosen sufficiently large. For example, if we take $\sigma=1$ and $n=0$, then we obtain the convergence speed $\Osym(k^{-1/2})$, which is not satisfactory for numerical calculations. In a next step, we need to solve the nonlinear system of equations $\Theta(\lambda,\mu)=0$ and $\hat\Theta(\lambda,\mu)=0$ with an appropriate numerical method, and then we get the eigenvalues $(\lambda,\mu)$ of the ellipsoidal wave equation \eqref{AEWE}. It is a further benefit of the approach described here that with the series coefficients $d_k$, $\tilde d_k$ and $\hat d_k$ generated to calculate $\Theta$ and $\hat\Theta_k$, we immediately obtain the corresponding ellipsoidal wave functions. Due to \eqref{Shift}, \eqref{HolSol0}, \eqref{FloqSol} and $w_{0,\rho}(z)=\Omega w_{1,\sigma}(z)$, any eigenfunction $w(z)$ of \eqref{AEWE} takes the form $w(z)=C_0 w_{0,\rho}(z)$ on $\mathfrak{D}_0$ and $w(z)=C_1 w_{1,\sigma}(z)$ on $\mathfrak{D}_1$, where
\begin{align*} 
w_{0,\rho}(z) & = z^{\alpha_0}(1-z)^{\beta_1+1}\spr{\eta_0(z)}{e_2}
= z^{-\rho/2}(1-z)^{(1+\sigma)/2}\sum_{k=0}^\infty (\spr{d_k}{e_2})z^k \\
w_{1,\sigma}(z) & = z^{\alpha_0}(1-z)^{\beta_2}\spr{\eta_2(z)}{e_2}
= z^{-\rho/2}(1-z)^{-\sigma/2}\sum_{k=0}^\infty (\spr{\tilde d_k}{e_2})(1-z)^k
\end{align*}
Moreover, according to \eqref{Trafo} and $\hat\Theta=0$, we have $w(z)=C_2 w_{2,\tau}(z)$ on $\mathfrak{D}_2$ with
\begin{equation*}
w_{2,\tau}(z)   
= \spr{\begin{pmatrix} -1 & 0 \\[1ex] 0 & 1 \end{pmatrix}\hat y_{0,\tau}(\hat z)}{e_2} 
= \left(\frac{c-z}{c-1}\right)^{-\tau/2}\left(\frac{z-1}{c-1}\right)^{(1+\sigma)/2}\sum_{k=0}^\infty (\spr{\hat d_k}{e_2})\left(\frac{c-z}{c-1}\right)^k
\end{equation*}
The constant factors have to be determined such that $C_0 w_{0,\rho}(z)=C_1 w_{1,\sigma}(z)$ holds on $\mathfrak{D}_0\cap\mathfrak{D}_1$ and $C_2 w_{2,\tau}(z)=C_1 w_{1,\sigma}(z)$ holds on $\mathfrak{D}_1\cap\mathfrak{D}_2$, where we can add an appropriate normalization of the eigenfunctions as a third condition.

We now consider a few examples with different parameters $\gamma$, $c$, $\rho$, $\sigma$, $\tau$ and pairs of eigenvalues $(\lambda,\mu)$, which we assume to be real, so that we may take $\Theta$ and also the eigenfunctions on $(0,1)$ and $(1,c)$ to be real functions.

As a first numerical example, we approximate $\Theta$ by $\Theta_k$ in the case $\gamma = 4$, $c = 1.6$, $\rho = 1$, $\sigma = 0$ (i.e. $\delta=0.5$) for the parameter values $\lambda = 3.2$, $\mu = -5$ with various $n$ and sufficiently large $k$. The calculation is terminated at index $k$ when the a posteriori estimate satisfied $\frac{k}{n+1.5}\,|\Theta_k-\Theta_{k-1}|\leq 10^{-10}$ for the first time, and this also provides an estimate for the error $|\Theta_k-\Theta|$. The results are listed in Table \ref{tab:Theta}. In the last column, the value $\Theta(3.2,-5) = -0.262836009163167617\ldots$ was calculated in advance with high precision, and it is used to compare the true error with the a posteriori estimate in the fourth column. The results in this example confirm that the number $k$ of computational steps that are required to meet the desired accuracy can be significantly reduced by slightly increasing $n$.

\begin{table}[htb]
\caption{Approximation of $\Theta$ by $\Theta_k$ using Theorem \ref{thm:AEWF} with $0\leq n\leq 5$, where the parameter values are $\lambda=3.2$, $\mu=-5$, $\gamma = 4$, $c = 1.6$, $\rho = 1$, $\sigma = 0$. All numerical values have been rounded to $18$ decimal places.} \label{tab:Theta}
\begin{tabular}{c||r|c|c|c} \toprule
$n$ & \multicolumn{1}{c|}{$k$} & $\Theta_k$ & $\frac{k}{n+1.5}\,|\Theta_k-\Theta_{k-1}|$ & $|\Theta_k-\Theta|$ \\ \midrule
 $0$ & $14510976$ & $-0.262836009063167646$ & $0.000000000099999993$ & $0.000000000099999971$ \\
 $1$ & $   28599$ & $-0.262836009063175185$ & $0.000000000099997249$ & $0.000000000099992432$ \\
 $2$ & $    1839$ & $-0.262836009063290010$ & $0.000000000099943779$ & $0.000000000099877607$ \\
 $3$ & $     358$ & $-0.262836009061005041$ & $0.000000000099207058$ & $0.000000000102162577$ \\
 $4$ & $     222$ & $-0.262836009256167290$ & $0.000000000097955235$ & $0.000000000092999672$ \\
 $5$ & $     154$ & $-0.262836009254764788$ & $0.000000000097901869$ & $0.000000000091597170$ \\
 \bottomrule
\end{tabular}
\end{table}

\begin{figure}[!htb]
\centering
\includegraphics{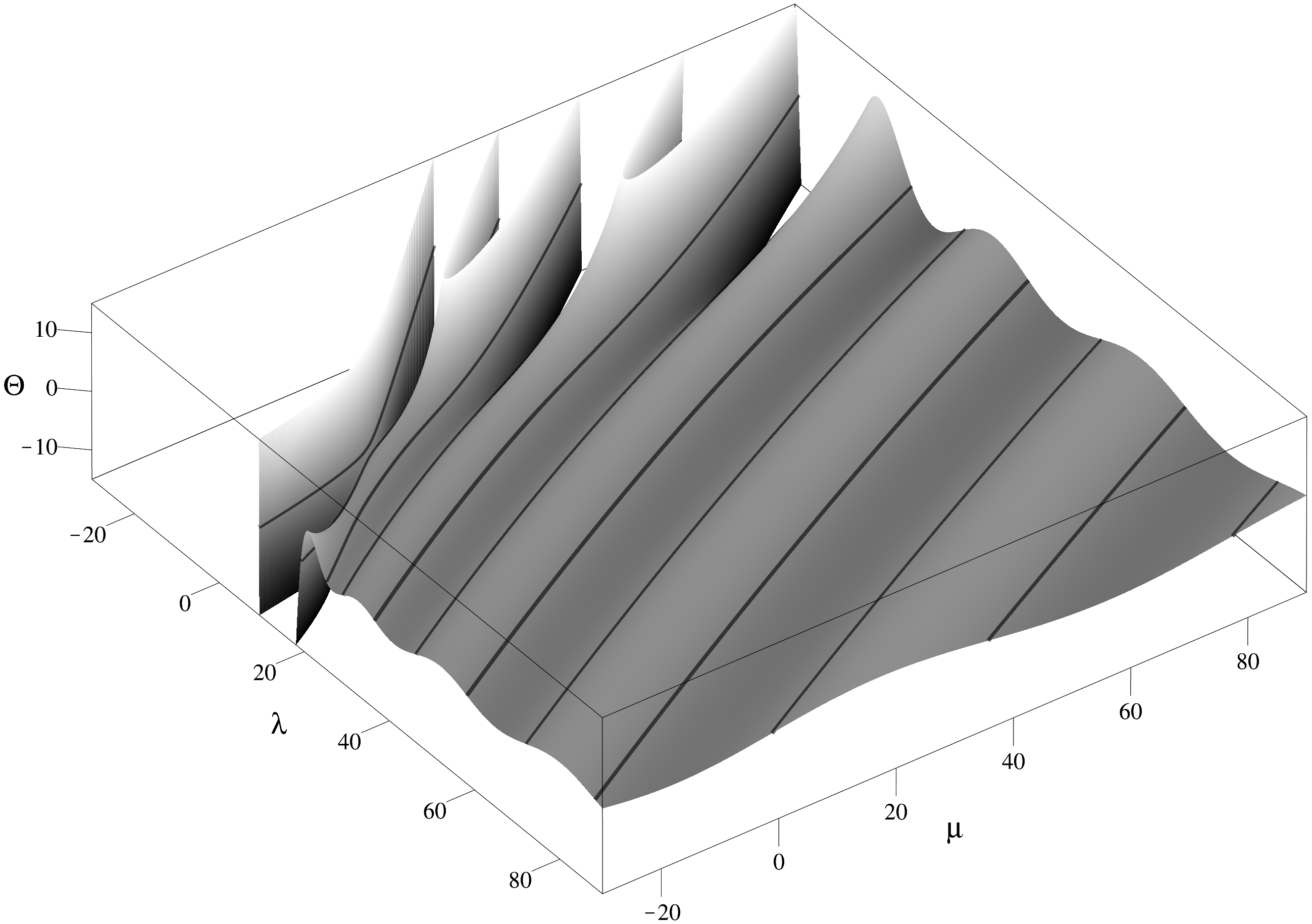}
\caption{The connection coefficient $\Theta$ for fixed $\gamma=4$, $c=1.6$, $\rho=1$, $\sigma=0$ as a function of the eigenvalue parameters $(\lambda,\mu)\in[-30,90]\times[-30,90]$. The contour lines indicate the zeros of the function $\Theta(\lambda,\mu)$.}
\label{fig:Theta}
\end{figure}

\begin{figure}[!htb]
\centering
\includegraphics{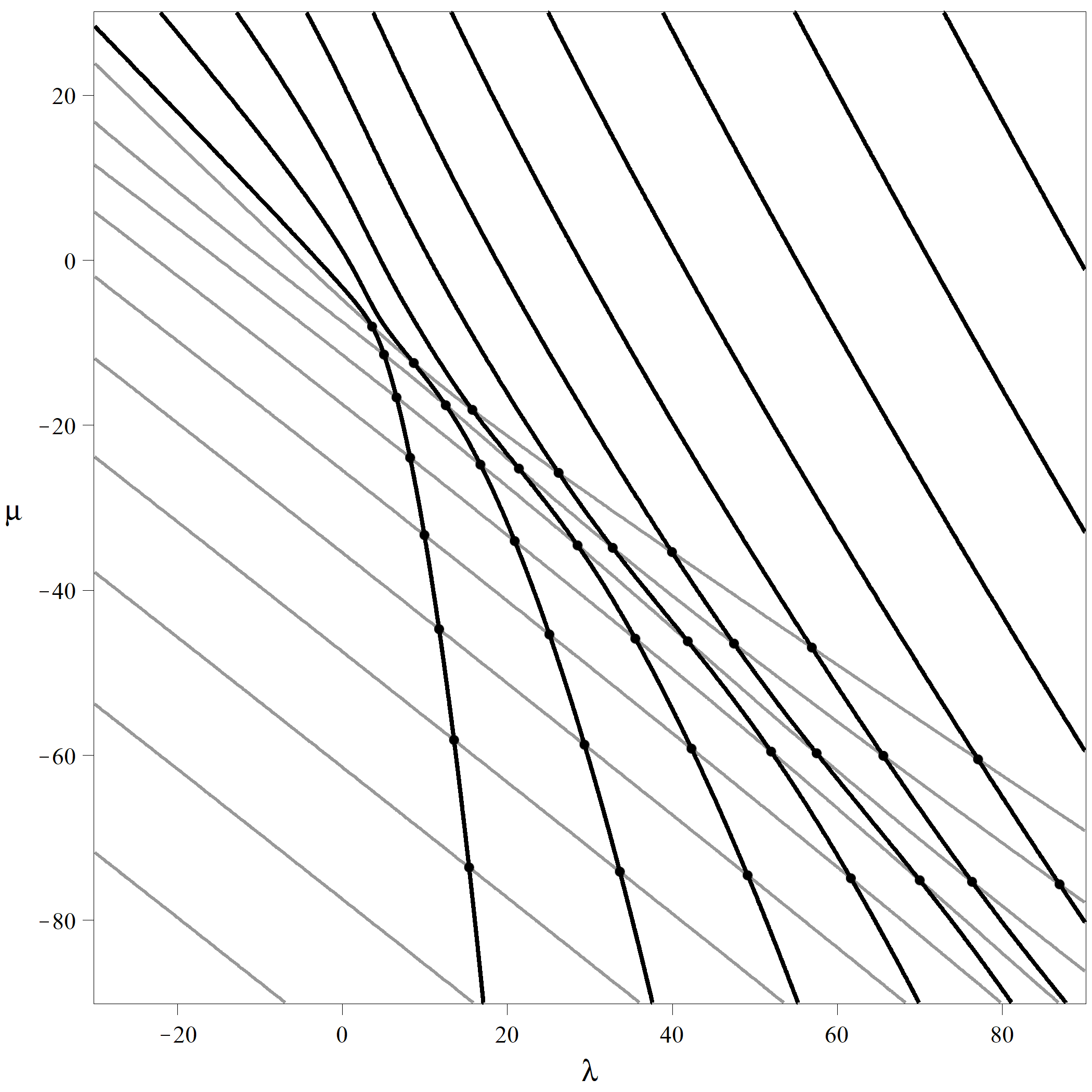}
\caption{The zeros of the connection coefficients $\Theta$ (black) and $\hat\Theta$ (gray) for the case $\gamma=4$, $c=1.6$, $\rho=1$, $\sigma=0$, $\tau=1$ in the domain $(\lambda,\mu)\in[-30,90]\times[-90,30]$. The intersection points of these curves, which are marked by small circles, correspond to the eigenvalues of the ellipsoidal wave equation for the parameters given above.}
\label{fig:EV}
\end{figure}

\begin{figure}[!htb]
\centering
\includegraphics[scale=0.8]{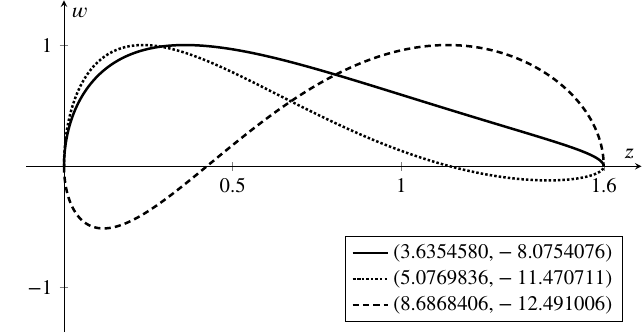} \\[1ex]
\includegraphics[scale=0.8]{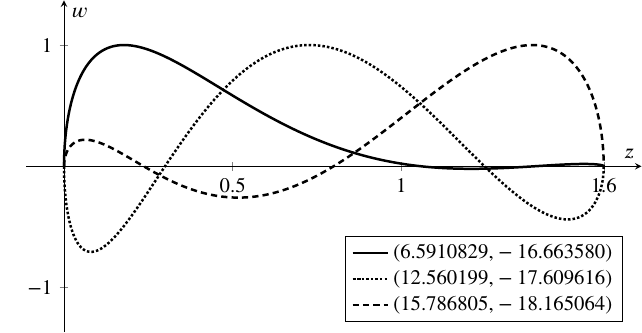}
\caption{The ellipsoidal wave functions in the case $\gamma=4$, $c=1.6$ and $\rho=1$, $\sigma=0$, $\tau=1$ on the interval $(0,c)$ with none, one and two zeros, each normalized to a maximum magnitude $1$. In the legend, the corresponding pairs of eigenvalues $(\lambda,\mu)$ are listed with an accuracy of eight decimal places.}\label{fig:AEWF}
\end{figure}

Figure \ref{fig:Theta} illustrates the function $\Theta=\Theta(\lambda,\mu)$ for the eigenvalue parameters in the domain $-30\leq\lambda\leq 90$ and $-30\leq\mu\leq 90$, again assuming $\gamma = 4$, $c = 1.6$, $\rho = 1$, $\sigma = 0$. In addition, the contour lines $\Theta=0$ are plotted there. For the parameters $(\lambda,\mu)$ along these zero curves, we obtain $w_{0,\rho} = \Omega w_{1,\sigma}$, i.e., these values satisfy the spectral problem on $\mathfrak{D}_0\cap\mathfrak{D}_1$. The contour plot in Figure \ref{fig:EV} not only shows the zeros of $\Theta(\lambda,\mu)$ but also the curves with $\hat\Theta(\lambda,\mu)=0$ in the case $\tau=1$. The eigenvalues of the ellipsoidal wave equation for $\gamma = 4$, $c = 1.6$, $\rho = 1$, $\sigma = 0$, $\tau=1$ are located at the intersection points of these curves. For the numerical calculation of these pairs of eigenvalues, the nonlinear system of equations $\Theta(\lambda,\mu)=0$ and $\hat\Theta(\lambda,\mu)=0$ were solved by means of a quasi-Newton method. To some of these eigenvalues the eigenfunctions are drawn in Figure \ref{fig:AEWF}. The ellipsoidal wave functions given there have been normalized so that their maximum magnitude is equal to $1$.

In Table \ref{tab:Lame}, some eigenvalues $(\lambda_n^m,\mu_n^m)$ of the Lam\'{e} equation (i.e., $\gamma=0$) have been calculated for the case $c=12/7$ and all eight possible combinations of $\rho$, $\sigma$, $\tau$. Here, $m$ and $n-m$ count the number of zeros of the corresponding Lam\'{e} function in $(0,1)$ and $(1,c)$, respectively. The number $c=12/7$ has been chosen to compare the eigenvalues with those given in \cite[Table 1]{WL:2005}. The results coincide, except for some cases (e.g., $\rho=0$, $\sigma=1$, $\tau=0$), where $\lambda_0^0$ in \cite{WL:2005} differs by a factor of $10$ from the value in given in Table \ref{tab:Lame}, but this seems to be just a typo in \cite{WL:2005}.

\begin{table}[htb]
\caption{Eigenvalues $(\lambda_n^m,\mu_n^m)$ of the Lam\'{e} equation ($\gamma=0$) for $c = \frac{12}{7}$ and the eight possible cases for $\rho$, $\sigma$, $\tau$.} \label{tab:Lame}
\begin{tabular}{c|c|c||r|r||r|r||r|r} \toprule \centering
$\rho$ & $\sigma$ & $\tau$ & $\lambda_0^0$ & $\mu_0^0$ & $\lambda_1^0$ & $\mu_1^0$ & $\lambda_1^1$ & $\mu_1^1$ \\ \midrule
$0$ & $0$ & $0$ & $0.000000$ &  $0.0$ & $0.611407$ & $-1.5$ & $2.102879$ & $-1.5$ \\
$0$ & $0$ & $1$ & $0.250000$ & $-0.5$ & $0.964286$ & $-3.0$ & $3.250000$ & $-3.0$ \\
$0$ & $1$ & $0$ & $0.428571$ & $-0.5$ & $0.981471$ & $-3.0$ & $4.304243$ & $-3.0$ \\
$1$ & $0$ & $0$ & $0.678571$ & $-0.5$ & $2.423953$ & $-3.0$ & $4.361761$ & $-3.0$ \\
$0$ & $1$ & $1$ & $0.678571$ & $-1.5$ & $1.303037$ & $-5.0$ & $5.482677$ & $-5.0$ \\
$1$ & $0$ & $1$ & $1.428571$ & $-1.5$ & $3.488893$ & $-5.0$ & $5.796821$ & $-5.0$ \\
$1$ & $1$ & $0$ & $1.964286$ & $-1.5$ & $3.597906$ & $-5.0$ & $7.473523$ & $-5.0$ \\
$1$ & $1$ & $1$ & $2.714286$ & $-3.0$ & $4.548506$ & $-7.5$ & $9.022923$ & $-7.5$ \\
\bottomrule
\end{tabular}
\end{table}

\begin{figure}[!htb]
\centering
\includegraphics{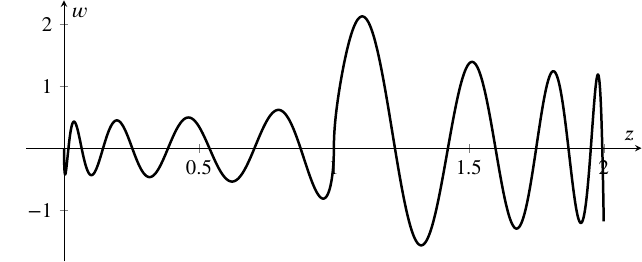}
\caption{The eigenfunction of the Lam\'{e} wave equation \eqref{Lame} for the parameters $(\rho,\sigma,\tau)=(1,1,0)$, $\omega^2=100$, $k^2=0.5$ and the eigenvalues $H=599.43708$, $L=629.53546$ satisfying the normalization condition \eqref{Norm}.}\label{fig:Fig-A8}
\end{figure}

\begin{figure}[!htb]
\centering
\includegraphics{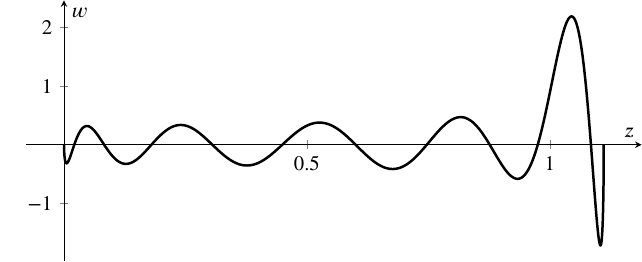}
\caption{The normalized eigenfunction of \eqref{Lame} satisfying \eqref{Norm} for the parameters $(\rho,\sigma,\tau)=(1,0,1)$, $\omega^2=1$, $k^2=0.9$ and the eigenvalues $(H,L)=(465.05152, 456.80932)$.}\label{fig:Fig-A9}
\end{figure}

As already mentioned, there are only a few publications that provide numerical data on the ellipsoidal wave equation, and one of them is \cite{ADKL:1989}. Furthermore, there is also no commonly accepted notation for its parameters. In \cite{Fedoryuk:1989a} and \cite{ADKL:1989}, the Lam\'{e} wave equation is written in the form
\begin{equation} \label{Lame}
z(z-1)(z-k^{-2})w''(z) + \tfrac{1}{2}\big(3z^2-2(1+k^{-2})z+k^{-2}\big)w'(z)+\tfrac{1}{4}\big(Hk^{-2}-Lk^{-2}z+\omega^2 z^2\big)w(z)=0
\end{equation} 
with spectral parameters $(H,L)$ instead of $(\lambda,\mu)$. Comparing \eqref{Lame} with our notation \eqref{AEWE} yields
\begin{equation} \label{HL}
c=k^{-2},\quad\gamma=\tfrac{1}{4}\omega^2,\quad\lambda=\tfrac{1}{4}Hk^{-2},\quad\mu=-\tfrac{1}{4}Lk^{-2}
\end{equation}
Table \ref{tab:ADKL} lists some eigenvalues for the case $\rho=1$, $\sigma=0$, $\tau=1$ at various $k^2$ and $\omega^2$. The numbers $(h,l)$ were taken from \cite[Table 1]{ADKL:1989}, while the eigenvalues $(H,L)$ have been calculated to four decimal places using our method (i.e., by solving $\Theta=\hat\Theta=0$) for the parameters given in \eqref{HL}. The pairs $(H,L)$ are in a good accordance with the data $(h,l)$ from \cite{ADKL:1989} except for minor deviations in the last decimal places. However, our results $(H,L)$ ought to be more reliable, since for these the values $\Theta$ and $\hat\Theta$ are closer to zero than for $(h,l)$. Figures \ref{fig:Fig-A8} and \ref{fig:Fig-A9} illustrate some eigenfunctions of \eqref{Lame} for the parameters and eigenvalues given in the captions. The eigenfunctions have been normalized such that
\begin{equation} \label{Norm}
\int_0^1\int_1^c\frac{(y-x) w(x)^2 w(y)^2}{|\varphi(x)|^{1/2}|\varphi(y)|^{1/2}}\,\mathrm{d}y\,\mathrm{d}x = 1,\quad\varphi(z) := z(1-z)(c-z)
\end{equation}
holds, as suggested in \cite[Eq. (2.3)]{ADKL:1989}. The graphs are consistent with those shown in \cite[Figures 8 and 9]{ADKL:1989}.

\begin{table}[htb]
\caption{Some eigenvalues of the Lam\'{e} wave equation \eqref{Lame} for the case $\rho=1$, $\sigma=0$, $\tau=1$ in the notation used by Abramov et\,al. The pairs $(h,l)$ were taken from \cite{ADKL:1989}, whereas $H$ and $L$ have been calculated with the method presented in Section 3.} \label{tab:ADKL}
\begin{tabular}{c|c||c|c||c|c} \toprule
$k^2$ & $\omega^2$ & $h$ & $l$ & $H$ & $L$ \\ \midrule
$0.5$ & $ 1$ & $404.574$ & $254.150$ & $404.5725$ & $254.1495$ \\
$0.5$ & $25$ & $415.437$ & $281.727$ & $415.4354$ & $281.7278$ \\
$0.5$ & $25$ & $105.656$ & $274.258$ & $105.6530$ & $274.2514$ \\
$0.5$ & $ 1$ & $102.028$ & $253.849$ & $102.0318$ & $253.8504$ \\
$0.9$ & $25$ & $141.090$ & $482.513$ & $141.0901$ & $482.5134$ \\
$0.9$ & $ 1$ & $137.683$ & $456.487$ & $137.6824$ & $456.4856$ \\
$0.9$ & $ 1$ & $465.062$ & $456.820$ & $465.0515$ & $456.8093$ \\
$0.9$ & $25$ & $476.765$ & $490.671$ & $476.7548$ & $490.6641$ \\
\bottomrule
\end{tabular}
\end{table}

\section{Application to the spheroidal wave equation}

As a further example, let us consider the spheroidal wave equation
\begin{equation} \label{ASWE}
\D{}{x}\left((1-x^2)\D{}{x}w(x)\right) + \left(\lambda + \gamma^2(1-x^2) - \frac{\mu^2}{1-x^2}\right)w(x) = 0,\quad -1<x<1
\end{equation}
where $\mu,\gamma^2\in\C$ are some fixed numbers, and $\lambda\in\C$ is considered to be the eigenvalue parameter. Throughout this section we assume that either $\mu=0$ or $\re(\mu)>0$ holds. The numbers $\lambda\in\C$ for which \eqref{ASWE} has a nontrivial bounded solution $w(x)$ on $(-1,1)$ are called spheroidal eigenvalues, and the corresponding functions $w(x)$ are the spheroidal wave functions.

Due to their importance for physics and engineering (see \cite[Chapter 4]{MS:1954} or \cite{Wang:2017}, for example), over the years a variety of methods have been developed to calculate the eigenvalues of \eqref{ASWE}. A frequently used approach is based on a series expansion with associated Legendre functions, where the recursion relation for the coefficients results in a transcendental equation with continued fractions whose roots are the spheroidal eigenvalues. These roots can be approximated by the eigenvalues of a symmetric tridiagonal matrix, see e.g. \cite{Hodge:1970}. In \cite{ADKPP:1984} a method is established which is based on the solution of an auxiliary differential equation for the associated phase function; it is similar to the algorithm used for the ellipsoidal wave equation in \cite{ADKL:1989}. Another method proposed in \cite{Skoro:2015} aims at smoothly matching the Floquet solutions $(1-x^2)^{\mu/2}\sum_{k=0}^\infty s_k(1-x)^k$ and $(1-x^2)^{\mu/2}\sum_{k=0}^\infty u_k(1+x)^k$, each of which is bounded in the neighborhood of one of the singular points $x=1$ and $x=-1$, respectively. 

In a recent paper \cite{Schmid:2023}, a new approach has been suggested to compute the eigenvalues of a spheroidal wave equation. For each $t\in\C$, a sequence $\Theta_k=\Theta_k(t)$ with limit $\Theta(t)=\lim_{k\to\infty}\Theta_k(t)$ is constructed, giving rise to an entire function $\Theta:\C\longrightarrow\C$ whose zeros $t_n$ are related to the spheroidal eigenvalues via $\lambda_n=t_n+\mu(\mu+1)$. Here, $\Theta(t)$ is simply the connection coefficient between certain Floquet solutions of \eqref{ASWE} at the regular singular points $x=\pm 1$. From \cite[Theorem 2.1]{Schmid:2023} it follows that $\Theta_k = \Theta(t) + \Osym(k^{\varepsilon-\re(\mu)-2})$ as $k\to\infty$ with arbitrary small $\varepsilon>0$. In the following, we will refine this approach by introducing an accelerated sequence which converges faster to $\Theta(t)$. More precisely, we generate a sequence $\Theta_k$ satisfying $\Theta_k = \Theta(t) + \Osym(k^{-\re(\mu)-n-2})$ for any given nonnegative integer $n$. Because the spheroidal wave functions of order zero play an important role in practical applications (see \cite{KM:2008} or \cite{ORX:2013}), we consider the special case $\mu=0$ in more detail and finally apply our results to some numerical examples. 

Initially, we associate a first order system to the second order ODE \eqref{ASWE}. To this end, we introduce
\begin{equation} \label{AVF}
y(z) = \begin{pmatrix} 2w'(2z-1) + \frac{\mu(2z-1)}{2z(1-z)}\,w(2z-1) \\[1ex] w(2z-1) \end{pmatrix},
\quad z\in(0,1)
\end{equation}
and the shifted eigenvalue parameter $t := \lambda-\mu(\mu+1)$. In \cite{Schmid:2023} it is shown that
$w(x)$ is a bounded solution of \eqref{ASWE} if and only if the associated vector function \eqref{AVF} is a bounded solution of the linear differential system
\begin{equation} \label{SphSys}
y'(z) = \left(\frac{1}{z}\begin{pmatrix} -\frac{\mu}{2}-1 & -t \\[1ex] 0 & \frac{\mu}{2} \end{pmatrix} + \frac{1}{z-1}\begin{pmatrix} -\frac{\mu}{2}-1 & t \\[1ex] 0 & \frac{\mu}{2} \end{pmatrix} + \begin{pmatrix} 0 & -4\gamma^2 \\[1ex] 1 & 0 \end{pmatrix}\right)y(z)
\end{equation}
on the interval $(0,1)$. This system has the form \eqref{GenSys}, where
\begin{equation*}
A := \begin{pmatrix} -\frac{\mu}{2}-1 & -t \\[1ex] 0 & \frac{\mu}{2} \end{pmatrix},\quad
B := \begin{pmatrix} -\frac{\mu}{2}-1 & t \\[1ex] 0 & \frac{\mu}{2} \end{pmatrix},\quad
G(z) \equiv G_0 := \begin{pmatrix} 0 & -4\gamma^2 \\[1ex] 1 & 0 \end{pmatrix}
\end{equation*}
Furthermore, if we define
\begin{equation*}
a_0 := \begin{pmatrix} \frac{-t}{\mu+1} \\[1ex] 1 \end{pmatrix}\quad\mbox{and}\quad
K := \begin{pmatrix} -1 & 0 \\[1ex] 0 & 1 \end{pmatrix}
\end{equation*}
then $a_0$ is an eigenvector of $A$ for the eigenvalue $\alpha_0 := \mu/2$. In addition, from $KAK=B$, $KBK=A$, $KG_0K=-G_0$ and $K^{-1}=K$ it follows that $y(z)$ is a solution of $\eqref{SphSys}$ if and only if $Ky(1-z)$ solves \eqref{SphSys}. According to \cite[Lemma 3.1]{Schmid:2023} with $\sigma=-t$ and $\rho=1$, this system has a fundamental matrix
\begin{equation*}
Y_0(z) = \begin{pmatrix} \frac{1}{z} & 0 \\[1ex] 0 & 1 \end{pmatrix}H_0(z)
\begin{pmatrix} z^{-\mu/2} & 0 \\[1ex] 0 & z^{\mu/2} \end{pmatrix}
\begin{pmatrix} 1 & 0 \\[1ex] q\log z & 1 \end{pmatrix},\quad
H_0(0) = \begin{pmatrix} 1 & 0 \\[1ex] p & 1 \end{pmatrix}
\end{equation*}
on the unit disk $\mathfrak{D}_0:=\{z\in\C:|z|<1\}$ with some holomorphic matrix function $H_0:\mathfrak{D}_0\longrightarrow\MzC$, where $q=0$ for the case $\mu\not\in\Z$ and $p=-1/\mu$ for the case $\mu\neq 0$; if $\mu=0$, then $p=0$ and $q=1$. The fundamental matrix $Y_0(z)$ provides a Floquet solution $y_0(z) := Y_0(z)e_2 = z^{\mu/2}h_0(z)$ of \eqref{SphSys}, where $h_0:\mathfrak{D}_0\longrightarrow\Cz$ is a holomorphic vector function with $h_0(0) = a_0$. Another fundamental matrix for \eqref{SphSys} is
\begin{equation*}
Y(z) := KY_0(1-z) = \begin{pmatrix} \frac{1}{1-z} & 0 \\[1ex] 0 & 1 \end{pmatrix}KH_0(1-z)
\begin{pmatrix} (1-z)^{-\mu/2} & 0 \\[1ex] 0 & (1-z)^{\mu/2} \end{pmatrix}
\begin{pmatrix} 1 & 0 \\[1ex] q\log(1-z) & 1 \end{pmatrix}
\end{equation*}
which can be converted into the following form:
\begin{equation*}
Y(z) = H(z)
\begin{pmatrix} (1-z)^{-\mu/2-1} & 0 \\[1ex] 0 & (1-z)^{\mu/2} \end{pmatrix}
\begin{pmatrix} 1 & 0 \\[1ex] q\log(1-z) & 1 \end{pmatrix}
\end{equation*}
where 
\begin{equation*}
H(z) := 
\begin{pmatrix} \frac{1}{1-z} & 0 \\[1ex] 0 & 1 \end{pmatrix}KH_0(1-z)
\begin{pmatrix} 1-z & 0 \\[1ex] 0 & 1 \end{pmatrix}\quad\mbox{with}\quad
H(1) = \begin{pmatrix} -1 & \ast \\[1ex] 0 & 1 \end{pmatrix}
\end{equation*}
is a holomorphic matrix function on $\mathfrak{D}_1:=\{z\in\C:|z-1|<1\}$. Hence, we can write $y_0(z)$ as a linear combination
\begin{equation*}
y_0(z) = Y(z)c,\quad c = c(t) = \begin{pmatrix} \Theta(t) \\[0.5ex] \Omega(t) \end{pmatrix}
\end{equation*}
on $\mathfrak{D}_0\cap\mathfrak{D}_1$. Finally, \cite[Lemma 3.3]{Schmid:2023} implies that $\lambda=t+\mu(\mu+1)$ is an eigenvalue of the spheroidal wave equation \eqref{ASWE} if and only if $\Theta(t)=0$; in this case, $y_0(z)= z^{\mu/2}h_0(z)$ is a constant multiple of $y_2(z):=Y(z)e_2=(1-z)^{\mu/2}Kh_0(1-z)$.

Now that we have seen that the spheroidal eigenvalues coincide with the zeros of the entire function $\Theta:\C\longrightarrow\C$, all we need is an effective algorithm that makes the calculation of $\Theta$ as simple as possible. At first we note that $A-\alpha_0-k=A-\mu/2-k$ is invertible for all positive integers $k$. Moreover, $\beta_1 := -\mu/2-1$ and $\beta_2 := \mu/2$ are the eigenvalues of $B$ with corresponding eigenvectors $b_1 := -e_1$ and $b_2 := K a_0$, respectively, which are as well the column vectors of $H(1)$. Finally, $\delta = \beta_2-\beta_1 = \mu+1$ satisfies $\re(\delta)>-1$, and therefore the assumptions (a) -- (c) in Section 1 are fulfilled. The transformation \eqref{Shift} turns \eqref{SphSys} into the differential system \eqref{RegSing} with $A_0 := A-\mu/2$ and $A_1 := B+\mu/2$. It has a holomorphic solution $\eta_0(z)=(1-z)^{\mu/2}z^{-\mu/2}y_0(z)=\sum_{k=0}^\infty z^k d_k$ on $\mathfrak{D}_0$, where $d_0 = a_0$, and a Floquet solution $\eta_2(z)=(1-z)^{\mu/2}z^{-\mu/2}y_2(z)=(1-z)^{\mu}\sum_{k=0}^\infty (1-z)^k\tilde d_k$ on $\mathfrak{D}_1$ with $\tilde d_0=Ka_0$. Lemma \ref{thm:RecRel} provides an algorithm for determining the coefficients $d_k$ and $\tilde d_k$, which we can use to calculate the connection coefficient $\Theta$ according to Theorem 4. All in all we get the following result:

\begin{Theorem} \label{thm:SphCon}
Assume $\re(\mu)>0$ or $\mu=0$, and let $n$ be a fixed nonnegative integer. In addition, let $d_k\in\Cz$ be a sequence of vectors which are computed by the recurrence relation
\begin{equation} \label{Coeff0}
\begin{split}
u_k & := \begin{pmatrix} 0 & \frac{t(\mu+1-k)}{k(\mu+1+k)} \\[1ex] 0 & -\frac{\mu+1}{k} \end{pmatrix} d_{k-1} - \begin{pmatrix} \frac{t}{k(\mu+1+k)} & \frac{4\gamma^2}{\mu+1+k} \\[1ex] -\frac{1}{k} & 0 \end{pmatrix} u_{k-1} \\
d_k & := d_{k-1}+u_k\quad\mbox{for}\quad k=1,2,3,\ldots\quad\mbox{with}\quad
u_0 = d_0 := \begin{pmatrix} \frac{-t}{\mu+1} \\[1ex] 1 \end{pmatrix}
\end{split}\end{equation}
and let the vectors $\tilde d_0,\tilde d_1,\ldots,\tilde d_n$ be given by means of
\begin{equation*}
\begin{split}
\tilde u_k & := \begin{pmatrix} \frac{\mu}{\mu+1+k} & \frac{t(k-1)}{k(\mu+1+k)} \\[1ex] 0 & -\frac{1}{k} \end{pmatrix}\tilde d_{k-1} + \begin{pmatrix} \frac{-t}{k(\mu+1+k)} & \frac{4\gamma^2}{\mu+1+k} \\[1ex] -\frac{1}{k} & 0 \end{pmatrix}\tilde u_{k-1} \\
\tilde d_k & := \tilde d_{k-1}+\tilde u_k\quad\mbox{for}\quad k=1,2,\ldots,n\quad\mbox{with}\quad
\tilde u_0 = \tilde d_0 := \begin{pmatrix} \frac{t}{\mu+1} \\[1ex] 1 \end{pmatrix}
\end{split}\end{equation*}
Finally, we introduce the vectors
\begin{equation*}
p_k := \tilde d_0 + \sum_{\ell=1}^n \left(\prod_{m=0}^{\ell-1}\frac{m+\mu}{m+\mu-k}\right)\tilde d_\ell,\quad
\vartheta_k := \frac{1}{e_2\T J p_k}J p_k,\quad J := \begin{pmatrix} 0 & 1 \\[1ex] -1 & 0 \end{pmatrix}
\end{equation*}
where $b_1\T Jp_k\neq 0$ for sufficiently large $k\geq k_1$. If we define $\Theta_k := \spr{d_k}{\vartheta_k}$ for all $k\geq k_1$, then $\lim_{k\to\infty}\Theta_k = \Theta$. More precisely, we get
$\Theta_k = \Theta + \Osym(k^{-\re(\mu)-n-2})$ as $k\to\infty$, and
\begin{equation*}
|\Theta - \Theta_k| \leq (1+\varepsilon)\frac{k\,|\Theta_k-\Theta_{k-1}|}{\re(\mu)+n+2}
\end{equation*}
holds for arbitrary given $\varepsilon>0$ and sufficiently large $k\geq k_2$.
\end{Theorem}

Theorem \ref{thm:SphCon} has already been proved in \cite{Schmid:2023}, but only for the special case $n=0$. In the proof given there, the coefficients $d_k$ were multiplied from the left by $\vartheta\T$, where $\vartheta$ is a vector orthogonal to $b_2$ with first component $1$, which is independent of $k$. The asymptotic behavior \eqref{LimCoeff} for $n=0$ yields $\spr{d_k}{\vartheta} = \Theta + \Osym(k^{-\re(\mu)-2})$ as $k\to\infty$, and the slowest speed $\Osym(k^{-2})$ for the convergence of $\spr{d_k}{\vartheta}$ to $\Theta$ is obtained in the case $\mu=0$, that is, for the spheroidal wave equation of order zero, which plays an important role in technical applications. In signal processing, for example, the prolate spheroidal wave functions of order zero are the eigenfunctions of a time-limiting operation followed by a low-pass filter (see \cite{SP:1961}). Now, by increasing $n$ and replacing $\vartheta$ with a slightly more sophisticated vector $\vartheta_k$, we are able to enhance the speed of convergence also in the case $\mu=0$.

\begin{Remark}
If $\mu=0$, then the computation of $\Theta$ in Theorem \ref{thm:SphCon} can be simplified to some extent, as $\eta_2(z) = Kh_0(1-z)=\sum_{k=0}^\infty(1-z)^k Kd_k$ and therefore $\tilde d_k = K d_k$ for all $k$. So we do not need to calculate the coefficients $\tilde d_k$ of $\eta_2$ separately.
\end{Remark}

From our previous considerations it follows that the series coefficients $d_k$, $\tilde d_k$ can be used to obtain the connection coefficient $\Theta=\Theta(t)$ as well as the Floquet solutions itself. If $\Theta(t)=0$, then $\lambda=t+\mu(\mu+1)$ is an eigenvalue of \eqref{ASWE}, and $\eta_0$ provides the corresponding spheroidal wave function. In summary, we get the following result:

\begin{Corollary} \label{thm:ASWF}
For each $t\in\C$, the limit $\Theta(t)=\lim_{k\to\infty}\Theta_k(t)$ of the sequence \eqref{Theta} exists and gives rise to an entire function $\Theta:\C\longrightarrow\C$ whose zeros $t_1,t_2,t_3,\ldots$ are related to the spheroidal eigenvalues by $\lambda_n=t_n+\mu(\mu+1)$. Moreover, from the vectors $d_k$ calculated in \eqref{Coeff0} with $t=t_n$ we receive the spheroidal wave function
\begin{equation*}
w(x) 
= \left(\tfrac{1+x}{1-x}\right)^{\mu/2}\sum_{k=0}^\infty\tfrac{e_2\T d_k}{2^k}(1+x)^k
= \Omega\left(\tfrac{1-x}{1+x}\right)^{\mu/2}\sum_{k=0}^\infty\tfrac{e_2\T d_k}{2^k}(1-x)^k,\quad -1<x<1
\end{equation*}
Here $\Omega=1$ or $\Omega=-1$, and hence $w(x)$ is either an even or an odd function.
\end{Corollary}

\begin{Proof}
It remains to prove the series expansion for the spheroidal wave functions in the case $t=t_n$. From \eqref{AVF}, \eqref{Shift}, and \eqref{HolSol0} it follows that
\begin{equation*}
w(x) := e_2\T y_0\big(\tfrac{1+x}{2}\big) 
= e_2\T\left(\tfrac{1+x}{1-x}\right)^{\mu/2}\eta_0\big(\tfrac{1+x}{2}\big)
= \left(\tfrac{1+x}{1-x}\right)^{\mu/2}\sum_{k=0}^\infty\tfrac{e_2\T d_k}{2^k}(1+x)^k
\end{equation*}
is a nontrivial bounded solution on $(-1,1)$. Moreover, if $t=t_n$ is a zero of $\Theta$, then $\eta_0(z)=\Omega\,\eta_2(z)$ and hence also 
\begin{equation} \label{SymRel}
y_0(z) = \Omega y_2(z) = \Omega K Y_0(1-z)e_2
= \begin{pmatrix} -\Omega & 0 \\[1ex] 0 & \Omega \end{pmatrix}y_0(1-z)
\end{equation}
If we put $z=1/2$ and take into account that $y_0(1/2)$ is not the zero vector, then we get $\Omega=\pm 1$. In addition, \eqref{SymRel} with $z=\frac{1+x}{2}$ yields $w(x) = e_2\T y_0\big(\tfrac{1+x}{2}\big) = \Omega e_2\T y_0\big(\tfrac{1-x}{2}\big) = \Omega\,w(-x)$. Thus, the spheroidal wave function $w(x)$ is either an even or an odd function, depending on whether $\Omega=1$ or $\Omega=-1$ holds. Note that the splitting of the spheroidal wave functions into even and odd functions is well known and simply results from the fact that for a solution $w(x)$ of \eqref{ASWE} also $w(-x)$ solves this differential equation.
\end{Proof}

As a numerical example, we approximate $\Theta(1.5)$ by $\Theta_k$ for fixed $\mu=0$ and $\gamma = 2$ with different $n$. The computation is terminated as soon as $\frac{k}{n+2}\,|\Theta_k-\Theta_{k-1}|\leq 10^{-12}$ is met, so that the absolute error $|\Theta_k-\Theta(1.5)|$ is also about $10^{-12}$ (provided $k$ is sufficiently large). The results have been rounded to $18$ decimal places, and they are listed in Table \ref{tab:Theta0}, where in the last column the value $\Theta(1.5)=0.349852604826025926...$ has been calculated in advance with high precision to compare the true error with the a posteriori estimate in the fourth column. Just as for the ellipsoidal wave equation in Section 3, the results confirm that one can accelerate the calculation significantly by slightly increasing $n$.

\begin{figure}
\centering
\includegraphics{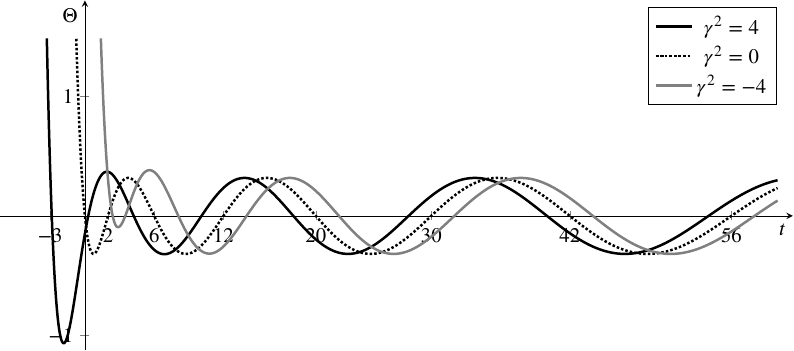}
\caption{The function $\Theta(t)=\Theta(t;\mu;\gamma^2)$ for $\mu=0$ and various $\gamma^2$. If $\gamma^2=0$, then \eqref{ASWE} reduces to the Legendre differential equation with eigenvalues $\lambda_n=n(n+1)$, $n\in\{0,1,2,3,\ldots\}$, and these are exactly the zeros of the dotted function $\Theta(t;0,0)$.}\label{fig:Theta0}
\end{figure}

Figure \ref{fig:Theta0} illustrates $\Theta(t)$ for $\mu=0$ and three different parameter values $\gamma^2$, where the graphs have been obtained by means of Corollary \ref{thm:ASWF}. Finally, if we calculate the zeros of the function $\Theta(t)$ in the case $\mu=0$ and $\gamma^2=4$ using, for example, the secant method, then we get the values listed in Table \ref{tab:EV0}, and these are also the eigenvalues of the corresponding spheroidal wave equation of order zero.

\begin{table}[htb]
\caption{Approximation of $\Theta(1.5)$ by $\Theta_k$ for $\mu=0$ in the prolate case $\gamma^2=4$ using Theorem \ref{thm:SphCon} with $0\leq n\leq 5$.} \label{tab:Theta0}
\begin{tabular}{c||r|c|c|c} \toprule
$n$ & \multicolumn{1}{c|}{$k$} & $\Theta_k$ & $\frac{k}{n+2}\,|\Theta_k-\Theta_{k-1}|$ & $|\Theta_k-\Theta(1.5)|$ \\ \midrule
 $0$ & $1928002$ & $0.349852604825025928$ & $0.000000000001000000$ & $0.000000000000999998$ \\
 $1$ & $   8133$ & $0.349852604825025664$ & $0.000000000000999731$ & $0.000000000001000262$ \\
 $2$ & $   2562$ & $0.349852604827023893$ & $0.000000000000999387$ & $0.000000000000997966$ \\
 $3$ & $    396$ & $0.349852604825045096$ & $0.000000000000999869$ & $0.000000000000980830$ \\
 $4$ & $    284$ & $0.349852604825050646$ & $0.000000000000995709$ & $0.000000000000975280$ \\
 $5$ & $     98$ & $0.349852604826806111$ & $0.000000000000922175$ & $0.000000000000780184$ \\
 \bottomrule
\end{tabular}
\end{table}

\begin{table}
\caption{The eigenvalues $\lambda_0,\ldots,\lambda_7$ of the spheroidal wave equation \eqref{ASWE} with $\mu=0$ and $\gamma^2=4$. They coincide with the numbers given in \cite[Table 1/1a]{SL:1964} for $c=\gamma=2.00$, taking into account that the eigenvalues are listed there in the form $\lambda_{0N} = \lambda_N+\gamma^2$.} \label{tab:EV0}
\begin{tabular}{c|r} \toprule
 $N$ & \multicolumn{1}{c}{$\lambda_N$} \\ \midrule
 $0$ & $-2.872265935150069$ \\
 $1$ & $ 0.287128543955796$ \\
 $2$ & $ 4.225713001105859$ \\
 $3$ & $10.100203876205334$ \\
 $4$ & $18.054829770465697$ \\
 $5$ & $28.035263096925295$ \\
 $6$ & $40.024747640293190$ \\
 $7$ & $54.018370784846266$ \\ \bottomrule
\end{tabular}
\end{table}

\section{Conclusion}

We have focused our examination of the linear differential system \eqref{GenSys} on the $2\times 2$ case, which plays an important role in applications and especially in the context of second-order ordinary differential equations. With a little more effort, one can generalize the results like the limit formula \eqref{LimCoeff} also to $n\times n$ systems. However, a more complicated distribution of the eigenvalues of $B$ and thus a more complex structure of the fundamental matrix $Y(z)$ has to be taken into account.

Furthermore, we can apply our results also to the generalized Heun equation
\begin{equation} \label{GHE}
w''(z) + \left(\frac{1-\nu_0}{z}+\frac{1-\nu_1}{z-1}+\frac{1-\nu_2}{z-c}+\kappa\right)w'(z)+
\frac{\lambda+\mu z+\gamma z^2}{z(z-1)(z-c)}w(z) = 0
\end{equation}
cf. \cite[Equation (0.3)]{SS:1980}, with fixed $c\in\C\setminus\{0,1\}$, $\gamma\in\C$ and some additional parameters $\kappa,\nu_1,\nu_2,\nu_3\in\C$ (the ellipsoidal wave equation is the special case where $c\in(1,\infty)$, $\nu_0=\nu_1=\nu_2=1/2$ and $\kappa=0$). For this purpose we just have to apply the transformation \eqref{Param} and replace \eqref{EllSys} by the differential system
\begin{equation*}
y'(z) = \left(\frac{1}{z}\begin{pmatrix} \nu_0-1 & a_{12} \\[1ex] 0 & 0 \end{pmatrix} 
+ \frac{1}{z-1}\begin{pmatrix} \nu_1-1 & b_{12} \\[1ex] 0 & 0 \end{pmatrix} 
+ \frac{1}{z-c}\begin{pmatrix} \nu_2-1 & r_{12} \\[1ex] 0 & 0 \end{pmatrix} 
- \frac{1}{c}\begin{pmatrix} \kappa c & 0 \\[1ex] 1 & 0 \end{pmatrix}\right)y(z)
\end{equation*}
which is again of type \eqref{GenSys}. Here, we may be interested in finding solutions of the form 
\begin{equation*}
w(z) = z^{\sigma_0}(z-1)^{\sigma_1}(z-c)^{\sigma_2}f(z)
\end{equation*}
for given $\sigma_j\in\{0,\nu_j\}$, $j\in\{0,1,2\}$, where $f:\C\longrightarrow\C$ is an entire function. If we assume $\nu_j\neq 1$ and $\re(\nu_j)>0$ for simplicity, then we can proceed as in Theorem \ref{thm:AEWF} and form sequences that converge to certain connection coefficients $\Theta$ and $\hat\Theta$, respectively, where the common zeros of $\Theta(\mu,\lambda)$ and $\hat\Theta(\mu,\lambda)$ provide the eigenvalue pairs $(\lambda,\mu)$ of \eqref{GHE}.


\begin{thebibliography}{10}
\expandafter\ifx\csname url\endcsname\relax
  \def\url#1{\texttt{#1}}\fi
\expandafter\ifx\csname urlprefix\endcsname\relax\def\urlprefix{URL }\fi
\expandafter\ifx\csname href\endcsname\relax
  \def\href#1#2{#2} \def\path#1{#1}\fi

\bibitem{SS:1980}
R.~Sch\"{a}fke, D.~Schmidt, The connection problem for general linear ordinary
  differential equations at two regular singular points with applications in
  the theory of special functions, SIAM Journal on Mathematical Analysis 11~(5)
  (1980) 848--862.
\newblock \href {https://doi.org/10.1137/0511076} {\path{doi:10.1137/0511076}}.

\bibitem{Schaefke:1980}
R.~Sch\"{a}fke, The connection problem for two neighboring regular singular
  points of general linear complex ordinary differential equations, SIAM
  Journal on Mathematical Analysis 11~(5) (1980) 863--875.
\newblock \href {https://doi.org/10.1137/0511077} {\path{doi:10.1137/0511077}}.

\bibitem{Schmid:2023}
H.~Schmid, Computation of the eigenvalues for the angular and {Coulomb}
  spheroidal wave equation, Applied Numerical Mathematics 185 (2023) 101--119.
\newblock \href {https://doi.org/10.1016/j.apnum.2022.11.018}
  {\path{doi:10.1016/j.apnum.2022.11.018}}.

\bibitem{Wasow:1976}
W.~Wasow, Asymptotic expansions for ordinary differential equations, 2nd
  Edition, Robert E. Krieger Publishing Co., Inc., Huntington, New York, 1976.

\bibitem{BSW:2005}
D.~Batic, H.~Schmid, M.~Winklmeier, On the eigenvalues of the
  {Chandrasekhar}-{Page} angular equation, Journal of Mathematical Physics
  46~(1) (2005) 012504.
\newblock \href {https://doi.org/10.1063/1.1818720}
  {\path{doi:10.1063/1.1818720}}.

\bibitem{MOS:1966}
W.~Magnus, F.~Oberhettinger, R.~P. Soni, Formulas and Theorems for the Special
  Functions of Mathematical Physics, 3rd Edition, Springer, Berlin -- Heidelberg, 1966.
\newblock \href {https://doi.org/10.1007/978-3-662-11761-3}
  {\path{doi:10.1007/978-3-662-11761-3}}.

\bibitem{ATZ:1983}
F.~M. Arscott, P.~J. Taylor, R.~V.~M. Zahar, On the numerical construction of
  ellipsoidal wave functions, Mathematics of Computation 40~(161) (1983)
  367--380.
\newblock \href {https://doi.org/10.2307/2007381} {\path{doi:10.2307/2007381}}.

\bibitem{Fedoryuk:1989a}
M.~V. Fedoryuk, {Lam\'{e}} wave functions, Mathematics of the USSR-Izvestiya
  33~(1) (1989) 179--200. \newline
\newblock \href {https://doi.org/10.1070/IM1989v033n01ABEH000819}
  {\path{doi:10.1070/IM1989v033n01ABEH000819}}.

\bibitem{Fedoryuk:1989b}
M.~V. Fedoryuk, The {Lam\'{e}} wave equation, Russian Mathematical Surveys
  44~(1) (1989) 153--180. \newline
\newblock \href {https://doi.org/10.1070/RM1989v044n01ABEH002009}
  {\path{doi:10.1070/RM1989v044n01ABEH002009}}.

\bibitem{WL:2005}
M.~Willatzen, L.~C. L.~Y. Voon, Numerical implementation of the ellipsoidal
  wave equation and application to ellipsoidal quantum dots, Computer Physics
  Communications 171~(1) (2005) 1--18.
\newblock \href {https://doi.org/10.1016/j.cpc.2005.04.006}
  {\path{doi:10.1016/j.cpc.2005.04.006}}.

\bibitem{ADKL:1989}
A.~A. Abramov, A.~L. Dyshko, N.~B. Konyukhova, T.~Levitina, Evaluation of
  {Lam\'{e}} angular wave functions by solving auxiliary differential
  equations, USSR Computational Mathematics and Mathematical Physics 29~(3)
  (1989) 119--131.
\newblock \href {https://doi.org/10.1016/0041-5553(89)90158-4}
  {\path{doi:10.1016/0041-5553(89)90158-4}}.

\bibitem{MS:1954}
J.~Meixner, F.~W. Sch\"{a}fke, Mathieusche Funktionen und
  Sph\"{a}roidfunktionen, Springer, New York, 1954, in German.
\newblock \href {https://doi.org/10.1007/978-3-662-00941-3}
  {\path{doi:10.1007/978-3-662-00941-3}}.

\bibitem{Wang:2017}
L.-L. Wang, A review of prolate spheroidal wave functions from the perspective
  of spectral methods, Journal of Mathematical Study 50~(2) (2017) 101--143.
\newblock \href {https://doi.org/10.4208/jms.v50n2.17.01}
  {\path{doi:10.4208/jms.v50n2.17.01}}.

\bibitem{Hodge:1970}
D.~B. Hodge, Eigenvalues and eigenfunctions of the spheroidal wave equation,
  Journal of Mathematical Physics 11~(8) (1970) 2308.
\newblock \href {https://doi.org/10.1063/1.1665398}
  {\path{doi:10.1063/1.1665398}}.

\bibitem{ADKPP:1984}
A.~A. Abramov, A.~L. Dyshko, N.~B. Konyukhova, T.~V. Pak, B.~S. Pariiskii,
  Evaluation of prolate spheroidal function by solving the corresponding
  differential equations, USSR Computational Mathematics and Mathematical
  Physics 24~(1) (1984) 1--11.
\newblock \href {https://doi.org/10.1016/0041-5553(84)90110-1}
  {\path{doi:10.1016/0041-5553(84)90110-1}}.

\bibitem{Skoro:2015}
S.~L. Skorokhodov, Evaluation of eigenvalues and eigenfunctions of {Coulomb}
  spheroidal wave equation, Matematicheskoe modelirovanie 27~(7) (2015)
  111--116, in Russian.

\bibitem{KM:2008}
A.~Karoui, T.~Moumni, New efficient methods of computing the prolate spheroidal
  wave functions and their corresponding eigenvalues, Applied and Computational
  Harmonic Analysis 24~(3) (2008) 269--289.
\newblock \href {https://doi.org/10.1016/j.acha.2007.06.004}
  {\path{doi:10.1016/j.acha.2007.06.004}}.

\bibitem{ORX:2013}
A.~Osipov, V.~Rokhlin, H.~Xiao, Prolate Spheroidal Wave Functions of Order Zero
  -- Mathematical Tools for Bandlimited Approximation, Springer, New York,
  2013.
\newblock \href {https://doi.org/10.1007/978-1-4614-8259-8}
  {\path{doi:10.1007/978-1-4614-8259-8}}.

\bibitem{SP:1961}
D.~Slepian, H.~O. Pollak, Prolate spheroidal wave functions, {Fourier} analysis
  and uncertainty -- {I}, Bell System Technical Journal 40~(1) (1961) 43--63.
\newblock \href {https://doi.org/10.1002/j.1538-7305.1961.tb03976.x}
  {\path{doi:10.1002/j.1538-7305.1961.tb03976.x}}.

\bibitem{SL:1964}
M.~M. Stuckey, L.~L. Layton, Numerical determination of spheroidal wave
  function eigenvalues and expansion coefficients, David Taylor Model Basin,
  Applied Mathematics Lab, Washington, D.C., 1964.

\end{thebibliography}

\end{document}